\documentclass{article}

\usepackage{arxiv}

\usepackage[utf8]{inputenc} 
\usepackage[T1]{fontenc}    
\usepackage{hyperref}       
\usepackage{url}            
\usepackage{booktabs}       
\usepackage{amsfonts}       
\usepackage{nicefrac}       
\usepackage{microtype}      
\usepackage{lipsum}
\usepackage{tikz}
\usetikzlibrary{arrows,shapes,quotes}
\usetikzlibrary{patterns,decorations.pathreplacing}
\usepackage{amssymb, graphicx}
\usepackage{hyperref,lscape}
\usepackage{blkarray}
\usepackage{color}
\usepackage{amsmath}
\usepackage{mathtools}
\usepackage{mathrsfs}
\usepackage{chngcntr}
\usepackage{apptools}
\usepackage{blkarray}
\usepackage{csquotes}
\DeclareMathOperator{\Ima}{Im}
\DeclareMathOperator\supp{supp}
\usepackage{parskip} 

\newtheorem{RE1}{Running example 1 - Part}
\newtheorem{RE2}{Running example 2 - Part}
\newtheorem{RE3}{Running example 3 - Part}

\newenvironment{Running example 1 - Part}[1][]{\begin{thm}[#1]\upshape}{\end{thm}}
\newenvironment{Running example 2 - Part}[1][]{\begin{thm}[#1]\upshape}{\end{thm}}

\usepackage{multicol}

\newtheorem{corollary}{Corollary}
\newtheorem{theorem}{Theorem}
\newtheorem{proposition}{Proposition}
\newtheorem{example}{example}
\newtheorem{definition}{Definition}
\newtheorem{remark}{Remark}
\newtheorem{lemma}{Lemma}

\title{Absolutely complex balanced kinetic systems}

\author{
  Editha C.~Jose \\
  Institute of Mathematical Sciences and Physics \\
  University of the Philippines \\ 
  Los Banos, Laguna 4031, Philippines \\
  \texttt{ecjose1@up.edu.ph} \\
   \And
Eduardo R.~Mendoza\thanks{Max Planck Institute of Biochemistry, Martinsried near Munich, Germany} \thanks{Faculty of Physics, Ludwig Maximilian University, Munich 80539 Germany} \\
  Mathematics and Statistics Department \\
  Center for Natural Sciences and Environmental Research \\
  De La Salle University \\
  Manila 0922, Philippines \\
  \texttt{eduardo.mendoza@dlsu.edu.ph} \\
   \And
  Dylan Antonio SJ.~Talabis \\
  Institute of Mathematical Sciences and Physics \\
  University of the Philippines \\
  Los Banos, Laguna 4031, Philippines \\
  \texttt{dstalabis1@up.edu.ph} \\

}

\begin{document}
\maketitle

\begin{abstract}
A complex balanced kinetic system is absolutely complex balanced (ACB) if every positive equilibrium is complex balanced. Two results on absolute complex balancing were foundational for modern chemical reaction network theory (CRNT): in 1972, M. Feinberg proved that any deficiency zero complex balanced system is absolutely complex balanced. In the same year, F. Horn and R. Jackson showed that the (full) converse of the result is not true: any complex balanced mass action system, regardless of its deficiency, is absolutely complex balanced. In this paper, we present initial results on the extension of the Horn and Jackson ACB Theorem. In particular, we focus on other kinetic systems with positive deficiency where complex balancing implies absolute complex balancing. While doing so, we found out that complex balanced power law reactant determined kinetic systems (PL-RDK) systems are not ACB. In our search for necessary and sufficient conditions for complex balanced systems to be absolutely complex balanced, we came across the so-called CLP systems (complex balanced systems with a desired "log parametrization" property). It is shown that complex balanced systems with bi-LP property are absolutely complex balanced. For non-CLP systems, we discuss novel methods for finding sufficient conditions for ACB in kinetic systems containing non-CLP systems: decompositions, the Positive Function Factor (PFF) and the Coset Intersection Count (CIC) and their application to poly-PL and Hill-type systems.
\end{abstract}

\keywords{chemical reaction network theory \and power law kinetics \and complex balancing \and absolute complex balancing }

\section{Introduction}
A positive equilibrium $x^*$ of a chemical kinetic system $(\mathscr{N}, K)$ is complex balanced if and only if $K(x^*)\in\ker I_a$, where $I_a$ is the incidence map of $\mathscr{N}.$  If the kinetic system is complex factorizable (CF), i.e., its species formation rate function $f(x) = Y\circ A_k \circ \Psi_K (x)$, where $A_k$ is a Laplacian, or equivalently $\Psi_K(x)\in\ker A_k.$ We denote the sets of positive and complex balanced equilibria of $(\mathscr{N}, K)$ with $E_+(\mathscr{N}, K)$ and $Z_+(\mathscr{N}, K)$ respectively. A kinetic system is complex balanced iff $Z_+(\mathscr{N}, K)\neq\emptyset.$ Since a network is weakly reversible if and only if $\ker I_a$ contains a positive element, the underlying network of a complex balanced system is necessarily weakly reversible. 

A complex balanced system is absolutely complex balanced (ACB) if every positive equilibrium is complex balanced, i.e., $E_+(\mathscr{N}, K)=Z_+(\mathscr{N}, K).$ Although the term was introduced only recently in \cite{FOME2020}, results on ``absolute complex balancing'' were foundational for Chemical Reaction Network Theory (CRNT) in the 1970's. In 1972, M. Feinberg showed that if a kinetic system with zero deficiency has a positive equilibrium, then the equilibrium is complex balanced. This implies that any deficiency zero complex balanced system is absolutely complex balanced (we refer to this result as the Feinberg ACB Theorem).  F. Horn, on the other hand, proved that any weakly reversible, deficiency zero mass action system is complex balanced, and hence also ACB. Furthermore, F. Horn and R. Jackson in the same year demonstrated that any complex balanced mass action system, regardless of deficiency, is absolutely complex balanced. We will call this result the Horn-Jackson ACB Theorem.

The Horn-Jackson ACB Theorem is interesting because it has at least three possible interpretations:
\begin{itemize}
    \item the ``Converse Counterexample'' (CC) view: it shows that the full converse to the Feinberg ACB Theorem is not true by providing a large set of ACB systems with deficiency $> 0$,
    \item the ``Partial Extension'' (PE) view: it extends the Feinberg ACB Theorem to positive deficiency systems but for a restricted set of kinetics (i.e. mass action only), and
    \item 	the ``Log Parametrization'' (LP) view: in his 1979 Wisconsin Lecture Notes, M. Feinberg showed that absolute complex balancing in a complex balanced mass action system is equivalent to its set of complex balanced equilibria being ``log parametrized'' by $S^\perp$, i.e., $Z_+(\mathscr{N}, K)  = \{ x \in \mathbb R^\mathscr S_> |\log x - \log x^*| \in S^\perp\}$, where $S$ is the network's stoichiometric subspace and $x^*$ a given complex balanced equilibrium. We call a complex balanced system with this ``log parametrization'' property a CLP system with flux space $S$ (and parameter space $S^\perp$).
\end{itemize}

The PE view leads to the following question (which we call the Horn-Jackson ACB Extension Problem or simply the Extension Problem): beyond mass action kinetics, which necessary or sufficient conditions ensure absolute complex balancing in systems with positive deficiency?  In particular, are there any kinetics sets other than mass action, where any complex balanced system with positive deficiency is absolutely complex balanced? We call this particular case the "Strong Extension Problem". This paper presents our initial results on the Extension Problem, which also reveal interesting (and partly surprising) connections to the CC and LP interpretations.

A natural candidate for extending the Horn-Jackson result is the complex balanced subset of power law kinetic systems with reactant-determined kinetic orders (denoted by PL-RDK), i.e., those where branching reactions of a reactant have identical rows in the system's kinetic order matrix. PL-RDK systems are precisely the complex factorizable power law systems and correspond to a subset of the generalized mass action systems (GMAS) introduced by S. M\"{u}ller and G. Rebensburger in 2014 \cite{MURE2014}.  A kinetic complex is the row in the kinetic order matrix of a reactant's reaction and, in analogy to the stoichiometric subspace, the kinetic order subspace $\tilde{S}$ is generated by the differences of the kinetic complexes of the complexes of reactions. In particular, M\"{u}ller and Regensburger showed that any complex balanced PL-RDK system is a CLP system with flux space  $=\tilde{S}$. Thus, the Extension Problem for PL-RDK is the same as whether the equivalence of ACB and CLP expressed in the LP view remains valid.

Our first important result is an example of a complex balanced PL-RDK system which is not ACB, providing a negative answer to the Strong Extension Problem for PL-RDK. While we provide a direct verification of $\emptyset\neq E_+(\mathscr{N}, K)\neq Z_+(\mathscr{N}, K)$, we note that considerations related to the CC view actually led us to the example (a detailed discussion is provided in Section \ref{sec:partialconverse}).

Building on results of M. Feinberg in his Lecture Notes, we broadened our study of the Extension Problem to the set of CLP systems. This approach also allowed us to address subsets of poly-PL and Hill-type systems with the CLP property studied in \cite{MHRM2020} and \cite{HEME2021} and have found applications, e.g., in evolutionary games with replicator dynamics \cite{TMMJ2020}. Our second important result is a necessary and sufficient condition, the bi-LP property, for ACB in CLP systems (Theorem  \ref{4.10}). This result provides a complete resolution of the Extension Problem in such systems, including all complex balanced PL-RDK systems. Clearly, connecting the PE and LP views is the basis for the result.

As initial steps in addressing the Extension Problem in non-CLP systems, we describe three methods for constructing ACB systems: (a) combinations of incidence independent and independent decompositions, (b) the Positive Function Factor (PFF) method and (c) the Coset Intersection Count (CIC) method. We illustrate with various kinetic systems how these methods can generate sufficient conditions for a complex balanced system to exhibit ACB. Furthermore, we derive a partial converse to the Feinberg ACB Theorem, i.e., identified the set of kinetic systems with kernel spanning equilibria images (KSE systems), which when absolutely complex balanced necessarily have zero deficiency. This result leads, on the one hand, to a necessary condition for ACB in a kinetic system with positive deficiency: it has to be non-KSE. On the other hand, it can be used to construct (if possible) a complex balanced system with positive deficiency which is not ACB. We illustrate the latter with how we originally found the example for PL-RDK systems.

In future work, we look to use these and other techniques to identify necessary and sufficient conditions for ACB in subsets of non-CLP systems to resolve the Extension Problem in those systems.

The paper is organized as follows: Section \ref{sec:fundamental} collects the fundamental concepts and results on chemical reaction networks and kinetic systems needed in the later sections. In Section \ref{sec:ACBPLRDK}, after the Extension Problem for the Horn-Jackson ACB Theorem is introduced and a counterexample for complex balanced PL-RDK systems is presented. Section \ref{sec:ACBCLP} broadens the scope to CLP systems and presents the necessary and sufficient condition resolving the Extension Problem for CLP systems. Sections \ref{sec:decompoACB} and \ref{sec:ACBHill} present the methods for finding sufficient conditions for ACB in kinetic systems containing non-CLP systems: decompositions, the Positive Function Factor (PFF) and the Coset Intersection Count (CIC) and their application to poly-PL and Hill-type systems. In Section \ref{sec:partialconverse}, a partial converse to the Feinberg ACB Theorem is derived and its use as a necessary condition and for constructing counterexamples for the Extension Problem is discussed. An overall summary is provided in Section \ref{sec:summary}.

\section{Fundamentals of chemical reaction networks and kinetic systems}\label{sec:fundamental}

In a chemical reaction, a $\textbf{species}$ is represented by a variable. We denote the nonempty finite set of distinct species by $\mathscr{S}=\left\{ X_1,X_2,...,X_m\right\}$ with cardinality of $\mathscr{S}$ equals $m$. A $\textbf{complex}$ is a linear combination of the species with nonnegative integer coefficients. We denote the nonempty finite set of complexes by $\mathscr{C}=\left\{ C_1,C_2,...,C_n\right\}$ where the cardinality of $\mathscr{C}$ is equal to $n$. A $\textbf{reaction}$ is an ordered pair of distinct complexes. Thus, if we denote this nonempty finite set of reactions by $\mathscr{R}$, we have  $\mathscr{R}\subset \mathscr{C} \times \mathscr{C}$. Let $r$ be the cardinality of $\mathscr{R}$. Consider the reaction
$$\alpha X_1 +\beta X_2 \rightarrow \gamma X_3,$$
$X_1$, $X_2$ and $X_3$ are the species. $\alpha X_1 +\beta X_2$ and $\gamma X_3$ are the complexes. In particular, $\alpha X_1 +\beta X_2$ is called the \textbf{reactant} (or \textbf{source}) \textbf{complex} and $\gamma X_3$ the \textbf{product complex}. The nonnegative coefficients $\alpha$, $\beta$ and $\gamma$ are called \textbf{stoichiometric coefficients}. Under power law kinetics (PLK), the rate at which the reaction occurs is given by
$$K = k X_1^a X_2^b$$
with rate constant $k > 0$ and $a,b \in \mathbb{R}$. We call $a$ and $b$ as kinetic orders. Thus, the reaction rate is a monomial in the reactant concentrations $X_1$ and $X_2$ with the exponents $a$ and $b$. Assuming mass action kinetics (MAK), (a subset of PLK), we have $a=\alpha$ and $b=\beta$, that is, the stoichiometric coefficients of the reactant complexes are the kinetic orders. PLK generalize MAK and has greater flexibility in modelling in biochemistry, epidemics, etc. \cite{SAV1998, ANMA1991}. Within a network involving additional species and reactions, the above reaction contributes to the dynamics of the species concentrations as
$$ \dot{X}=\left[ 
\begin{array}{c}
    \dot{X_1} \\
    \dot{X_2} \\
    \dot{X_3} \\
    \vdots \\
\end{array}
\right]=k X_1^a X_2^b \left( 
\begin{array}{c}
    -\alpha \\
    -\beta \\
    \gamma \\
    \vdots \\
\end{array}
\right) + \hdots$$

This is known as the \textbf{dynamical system} or system of ordinary differential equations (ODEs). See an example of PLK system below.

\begin{RE1}
\label{RE1a}
$$\begin{array}{l}
    \dot{X_1}= -k_1 X_1^{2} + k_2 X_1 X_2 - k_3 X_1 X_2 + k_4 X_3 - k_5 X_1^{2} + k_6 X_3^{2} - k_7 X_3^{2} + k_8 X_2^{-1} X_3^{-1} \\
    \dot{X_2}= k_1 X_1^{2} - k_2 X_1 X_2 + k_3 X_1 X_2 - k_4 X_3 \\
    \dot{X_3}= k_5 X_1^{2} - k_6 X_3^{2} + k_7 X_3^{2} - k_8 X_2^{-1} X_3^{-1}\\
\end{array}$$
where $k_i$'s are greater than $0$. 
\end{RE1}
Chemical reaction networks (CRNs) can be represented as a directed graph. The vertices or nodes are the complexes and the reactions are the edges. The CRN is not unique and might not have a physical interpretation. In Running example 1, the CRN can be the following: \\  

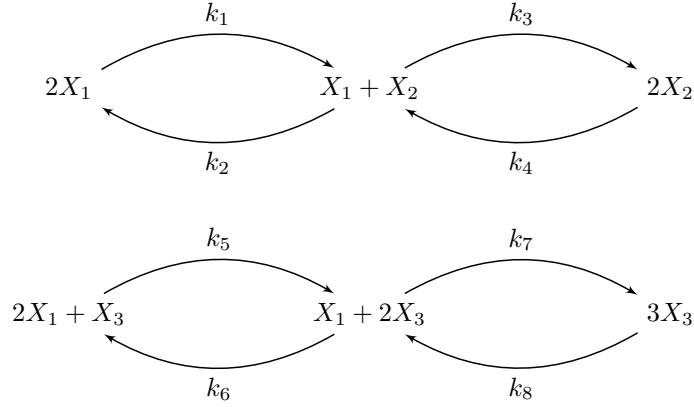
\begin{center}
\begin{tikzpicture}[baseline=(current  bounding  box.center)]
\tikzset{vertex/.style = {draw=none,fill=none}}
\tikzset{edge/.style = {bend left,->,> = latex', line width=0.20mm}} 
\node[vertex] (1) at  (-4.0,0) {$2X_1$};
\node[vertex] (2) at  (0.0,0) {$X_1 + X_2$};
\node[vertex] (3) at  (4.0,0) {$2X_2$};
\node[vertex] (4) at  (-4.0,-3) {$2X_1 + X_3$};
\node[vertex] (5) at  (0.0,-3) {$X_1 + 2X_3$};
\node[vertex] (6) at  (4.0,-3) {$3X_3$};
\draw [edge]  (1) to["$k_1$"] (2);
\draw [edge]  (2) to["$k_2$"] (1);
\draw [edge]  (2) to["$k_3$"] (3);
\draw [edge]  (3) to["$k_4$"] (2);
\draw [edge]  (4) to["$k_5$"] (5);
\draw [edge]  (5) to["$k_6$"] (4);
\draw [edge]  (5) to["$k_7$"] (6);
\draw [edge]  (6) to["$k_8$"] (5);
\end{tikzpicture}
\end{center}

The $k_i$'s are called the reaction rate constants. We have $m=3$ (species), $n=6$ (complexes), $n_r=6$ (reactant complexes) and $r=8$ (reactions). We can write
$$\mathscr{S}=\left\{ X_1, X_2, X_3\right\}, \quad \mathscr{C}=\left\{2X_1, X_1 + X_2, 2X_2, 2X_1+X_3, X_1+2X_3, 3X_3 \right\}.$$

On the other hand, the set of reaction $\mathscr{R}$ consists of the following:
$$\begin{array}{l}
R_{1}: 2X_1 \rightarrow X_1+X_2 \\
R_{2}: X_1+X_2 \rightarrow 2X_1\\
R_{3}: X_1+X_2 \rightarrow 2X_2\\
R_{4}: 2X_2 \rightarrow X_1+X_2\\
\end{array} \quad \begin{array}{l}
R_{5}: 2X_1+X_3 \rightarrow X_1+2X_3\\
R_{6}: X_1+2X_3 \rightarrow 2X_1+X_3\\
R_{7}: X_1+2X_3 \rightarrow 3X_3\\
R_{8}: 3X_3 \rightarrow X_1+2X_3\\
\end{array}$$

We denote the CRN $\mathscr{N}$ as $\mathscr{N}= (\mathscr{S}, \mathscr{C}, \mathscr{R})$. The \textbf{linkage classes} of a CRN are the subgraphs of a reaction graph where for any complexes $C_i$, $C_j$ of the subgraph, there is a path between them. Thus, the number of linkage classes, denoted as $l$, of Running example 1 is two ($l=2$). The linkage classes are:
$$\mathscr{L}_1=\left\{ R_1,R_2,R_3,R_4 \right\}, \quad \mathscr{L}_2=\left\{R_5,R_6,R_7,R_8 \right\}.$$
A subset of a linkage class where any two vertices are connected by a directed path in each direction is said to be a \textbf{strong linkage class}. Considering Running example \ref{RE1a}, there are two strong linkage classes whose number is denoted by $sl$. We also identify the \textbf{terminal strong linkage classes}, the number denoted as $t$, to be the strong linkage classes where there is no reaction from a complex in the strong linkage class to a complex outside the same strong linkage class.  The terminal strong linkage classes can be of two kinds: cycles (not necessarily simple) and singletons (which we call ``terminal points''). 

We now define important CRN classes.  A CRN is \textbf{weakly reversible} if every linkage class is a strong linkage class. A CRN is \textbf{t-minimal} if $t = l$, i.e. each linkage class has only one terminal strong linkage class. Let $n_r$ be the number of reactant complexes of a CRN. Then $n - n_r$ is the number of terminal points. A CRN is called \textbf{cycle-terminal} if and only if $n - n_r = 0$, i.e., each complex is a reactant complex. Clearly, the CRN of the Running example 1 is t-minimal weakly reversible. The dynamical system $f(x)$ (or species formation rate function (SFRF)) of the Running example 1 can be written as
$$\left[ 
\begin{array}{c}
    \dot{X_1} \\
    \dot{X_2} \\
    \dot{X_3} \\
\end{array}
\right]=
\begin{blockarray}{cccccccc}
R_1 & R_2 & R_3 & R_4 & R_5 & R_6 & R_7 & R_8  \\
\begin{block}{[cccccccc]}
-1 & 1 & -1 & 1 & -1 & 1 & -1 & 1 \\
1 & -1 & 1 & -1 & 0 & 0 & 0 & 0 \\
0 & 0 & 0 & 0 & 1 & -1 & 1 & -1 \\
\end{block}
\end{blockarray}
\left[ 
\begin{array}{c}
	k_1 X_1^{2} \\
	k_2 X_1 X_2 \\
    k_3 X_1 X_2 \\
	k_4 X_3 \\
	k_5 X_1^{2} \\
	k_6 X_3^{2} \\
	k_7 X_3^{2} \\
	k_8 X_2^{-1} X_3^{-1} \\
\end{array}
 \right]
=NK(x).$$
$N$ is called the stoichiometric matrix and $K(x)$ is called the kinetic vector. With each reaction $y\rightarrow y'$, we associate a \textbf{reaction vector} obtained by subtracting the reactant complex $y$ from the product complex $y'$. The \textbf{stoichiometric subspace} $S$ of a CRN is the linear subspace of $\mathbb{R}^\mathscr{S}$ defined by
$$S := \text{span } \left\lbrace y' - y \in \mathbb{R}^\mathscr{S} \mid y\rightarrow y' \in \mathscr{R}\right\rbrace.$$

The \textbf{map of complexes} $\displaystyle{Y: \mathbb{R}^\mathscr{C} \rightarrow \mathbb{R}^\mathscr{S}_{\geq}}$ maps the basis vector $\omega_y$ to the complex $ y \in \mathscr{C}$. 
The \textbf{incidence map} $\displaystyle{I_a : \mathbb{R}^\mathscr{R} \rightarrow \mathbb{R}^\mathscr{C}}$ is defined by mapping for each reaction $\displaystyle{R_i: y \rightarrow y' \in \mathscr{R}}$, the basis vector $\omega_{R_i}$ (or simply $\omega_i$) to the vector $\omega_{y'}-\omega_{y} \in \mathscr{C}$. 
The \textbf{stoichiometric map} $\displaystyle{N: \mathbb{R}^\mathscr{R} \rightarrow \mathbb{R}^\mathscr{S}}$ is defined as $N = Y \circ  I_a$. 

In Running example 1, the matrices $Y$ and $I_a$ are

$$Y=\begin{blockarray}{ccccccc}
C_1 & C_2 & C_3 & C_4 & C_5 & C_6 \\
\begin{block}{[cccccc]c}
2 & 1 & 0 & 2 & 1 & 0 & X_1 \\ 
0 & 1 & 2 & 0 & 0 & 0 & X_2 \\ 
0 & 0 & 0 & 1 & 2 & 3 & X_3 \\ 
\end{block}
\end{blockarray}$$

$$I_a=\begin{blockarray}{ccccccccc}
R_1 & R_2 & R_3 & R_4 & R_5 & R_6 & R_7 & R_8  \\
\begin{block}{[cccccccc]c}
-1 & 1 & 0 & 0 & 0 & 0 & 0 & 0 & C_1 \\
1 & -1 & -1 & 1 & 0 & 0 & 0 & 0 & C_2 \\
0 & 0 & 1 & -1 & 0 & 0 & 0 & 0 & C_3 \\
0 & 0 & 0 & 0 & -1 & 1 & 0 & 0 & C_4 \\
0 & 0 & 0 & 0 & 1 & -1 & -1 & 1 & C_5 \\
0 & 0 & 0 & 0 & 0 & 0 & 1 & -1 & C_6 \\
\end{block}
\end{blockarray}.$$
Here, we denote the complexes by $C_i, i=1,\dots,6$.

The \textbf{deficiency} $\delta$ is defined as $\delta = n - l - \dim S$. This non-negative integer is, as Shinar and Feinberg pointed out in \cite{SF2011}, essentially a measure of the linear dependency of the network's reactions. In Running example \ref{RE1a}, the deficiency of the network is 2. It is one of the important parameters in CRNT to establish claims regarding the existence, multiplicity, finiteness and parametrization of the \textbf{set of positive steady states}, denoted as $E_+$. It is defined as $E_+(\mathscr{N},K)=\left\{x\in \mathbb{R}^{\mathscr{S}}_> \middle| NK(x)=0 \right\}$. 

A CRN $\mathscr{N}$ with the property $\delta=\delta_1+\delta_2+\ldots+\delta_l$ is called a
\textbf{network with independent linkage classes} where $\delta_i$ is the deficiency of linkage class $\mathscr{L}_i$. Otherwise, if $\delta \neq \delta_1+\delta_2+\ldots+\delta_l$ it is called \textbf{dependent linkage class network}. A network is called conservative if $S^{\perp}$ contains a positive vector, where $S^{\perp}$ is the orthogonal complement of $S$. Horn and Jackson showed that a CRN is conservative if and only if $S$ is compact. In \cite{HOJA1972}, Horn and Jackson introduced a subset of $E_+$ called the set of complex balanced of equilibria denoted as $Z_+$. A kinetic system is complex balanced at a state (i.e. a species composition) if for each complex, formation and degradation are at equilibrium. A positive vector $c$ in $\mathbb{R}^{\mathscr{S}}$ is called \textbf{complex balanced} (CB) if $K(c)$ is contained in $\ker I_a$. A kinetic system is called \textbf{complex balanced} if it has a complex balanced equilibrium. The \textbf{set of complex balanced equilibria} of a kinetic system is defined as $Z_+(\mathscr{N},K)=\left\{x\in \mathbb{R}^{\mathscr{S}}_+ \middle| I_a K(x)=0 \right\}$. 

We introduce absolutely complex balanced kinetic system:

\begin{definition}
A kinetic system $(\mathscr{N},K)$ is \textbf{absolutely complex balanced (ACB)} if 
$$Z_+(\mathscr{N},K)\neq \emptyset \text{ and } E_+(\mathscr{N},K)=Z_+(\mathscr{N},K),$$
i.e., all positive equilibria are complex balanced.
\end{definition}

In this paper, we focus on CRNs endowed with \textbf{power law kinetics}. They have the form
$$\displaystyle K_{i}(x)=k_i \Pi^{r}_{i=1} x^{F_{ij}} \quad \forall i \in \left\{1,\ldots,r\right\}$$
with $k_i \in \mathbb{R}_+$ and $F_{ij} \in \mathbb{R}$. Power law kinetics is defined by an  $r \times m$ matrix $F=[F_{ij}]$, called the \textbf{kinetic order matrix}, and vector $k \in \mathbb{R}^r$, called the \textbf{rate vector}. In our Running example 1, the kinetic order matrix is
\begin{equation}
\nonumber
F=\begin{blockarray}{cccc}
X_1 & X_2 & X_3 \\
\begin{block}{[ccc]c}
2 & 0 & 0 & R_1 \\  
1 & 1 & 0 & R_2 \\  
1 & 1 & 0 & R_3 \\  
0 & 0 & 1 & R_4 \\ 
2 & 0 & 0 & R_5 \\  
0 & 0 & 2 & R_6 \\  
0 & 0 & 2 & R_7 \\  
0 & -1 & -1 & R_8 \\  
\end{block}
\end{blockarray}.
\end{equation}

A PLK system has \textbf{reactant-determined kinetics} (of type \textbf{PL-RDK}) if for any two reactions $R_i$, $R_j \in \mathscr{R}$ with identical reactant complexes, the corresponding rows of kinetic orders in $F$ are identical, i.e. $F_{ih}=F_{jh}$ for $h  \in \left\{1,\ldots,m \right\}$. A kinetics $K$ is \textbf{complex factorizable (CF)} if we can decompose $K=I_k \circ \Psi_K$ and $I_a \circ K =A_k \circ \Psi_K$ where $A_k$ is the Laplacian Map, $\Psi_K$ is the factor map, and $I_k$ is the interaction map. The kinetics of Running example 1 is complex factorizable with factor map $\Psi_K = \left[ \begin{array}{cccccc}
X_1^2 & X_1 X_2 & X_3 & X_1^2 & X_3^2 & X_2^{-1} X_3^{-1} \\ 
\end{array}
\right]^\top .$  

The concept of span surjectivity was introduced in \cite{AJLM2017} for any function $f: V \rightarrow W$ between finite dimensional real vector spaces. We say $f$ is span surjective $\Leftrightarrow$ $< \Ima f > = W$. The property is equivalent to the coordinate functions of $f$ being linearly independent (over $\mathbb{R}$). The property occurs in various situations, e.g. if $f$ is the "core" species formation rate function of a kinetics $K$, i.e. $f = NK : \Omega \rightarrow \mathscr{S}$, where $N$ is the stoichiometric matrix and $S = \Ima N$ is the stoichiometric subspace, then $f$ is span surjective $\Leftrightarrow$ the kinetic subspace of $K$ given by $< \Ima f >$ is equal to $S$. 

A complex factorizable (CF) kinetics $K$ is factor span surjective (FS) if its factor map $\Psi_K: \Omega \rightarrow \mathbb{R}^{\mathscr{C}}$ is span surjective.  M. Feinberg and F. Horn proved in \cite{FEHO1977} that mass action kinetics are factor span surjective. The property is often implicitly used in results about mass action kinetics. For example, in the important result on the coincidence of the kinetic and stoichiometric subspaces for mass action systems \cite{FEHO1977}, the main statement says: if a network is t-minimal then the kinetic subspace coincides with the stoichiometric subspace. As shown in \cite{AJLM2017}, this result extends only to t-minimal networks with FS kinetics, not to all CF kinetics.

M\"{u}ller and Regensburger \cite{MURE2014} 
introduced the $m \times n$ matrix $\widetilde{Y}$ for cycle terminal networks (that is, every complex is a reactant complex). This concept was generalized into arbitrary networks and was defined as follows.

\begin{equation}
\nonumber
(\widetilde{Y})_{ij}=\begin{cases}
  F_{ki}, \text{if } j \text{ is a reactant complex of reaction }  k, \\
  0, \text{otherwise} \\
  \end{cases}
\end{equation}
where $F$ is the kinetic order matrix. 

In 2018, Talabis et al. \cite{TAM2017} defined the T-matrix and the augmented T-matrix ($\widehat{T}$) as follows:

\begin{definition}
\label{def:25}
The $m \times n_r$ \textbf{T-matrix} is the truncated $\widetilde{Y}$ where the non-reactant columns are deleted and $n_r$ is the number of reactant complexes. The T-matrix defines a map $T: \mathbb{R}^{\rho(\mathscr{R})} \rightarrow \mathbb{R}^\mathscr{S}$. The \textbf{kinetic reactant subspace} $\widetilde{R}$ is the image of $T$. Its dimension is called the \textbf{kinetic reactant rank} $\widetilde{q}$.
\end{definition}

Define the $n_r \times l$ matrix $L = \left[ e_1,e_2,...,e_l \right]$ where $e^i$ is a characteristic vector for linkage class $\mathscr{L}^i$. The block matrix $\widehat{T} \in \mathbb{R}^{(m+l)\times n_{r}}$ is defined as 
\begin{equation}
\nonumber
\widehat{T}=\left[ 
\begin{array}[center]{c} T \\
L^{\top} \\
\end{array} \right].
\end{equation}

In \cite{TAM2017}, Talabis et al. defined the PL-TIK systems, a subclass of PL-RDK systems.

\begin{definition}
\label{TIK} 
A PL-RDK kinetics is  \textbf{$\widehat{T}$-rank maximal} (of type \textbf{PL-TIK}) if its column rank is maximal. 
\end{definition}

\begin{definition}
\label{kinetic reactant deficiency}
Let $\mathscr{N}$ be a network with $n_r$ reactant complexes and $K$ a PL-RDK kinetics with T-matrix $T$. If $\widehat{q}=\text{rank}(\widehat{T})$, then the \textbf{kinetic reactant deficiency} $\widehat{\delta}$ is defined as $$\widehat{\delta}= n_r - \widehat{q}.$$
\end{definition}

The kinetic reactant deficiency (a non-negative integer) measures the degree of the kinetic interactions of the PL-RDK system. The higher the kinetic reactant deficiency, the lower the extent of linear independence of kinetic orders (kinetic interaction). 

In our previous work \cite{TAME2018}, we proved a characterization of PL-TIK systems:
\begin{proposition}
PL-TIK systems and the zero kinetic reactant deficiency systems are equivalent. 
\end{proposition}

\begin{theorem}[\cite{TAME2018}]
\label{RKDZT}
Let $(\mathscr{N},K)$ be a PL-TIK system. $\mathscr{N}$ is weakly reversible if and only if $Z_+(\mathscr{N}, K) \neq \emptyset$.
\end{theorem}

\section{Complex balanced PL-RDK systems which are not absolutely complex balanced}\label{sec:ACBPLRDK}

Two results on absolute complex balancing in 1972 were foundational for Chemical Reaction Network Theory (CRNT). M. Feinberg showed that in a deficiency zero kinetic system with a positive equilibrium, every positive equilibrium is complex balanced. This implies the following:
\begin{theorem}[Feinberg ACB Theorem, \cite{FEIN1972a}]
Any complex balanced system $(\mathscr{N}, K)$ with zero deficiency is absolutely complex balanced.
\end{theorem}
F. Horn and R. Jackson obtained the following partial extension of the previous Theorem:
\begin{theorem}[Horn-Jackson ACB Theorem, \cite{HOJA1972}]\label{ACBmass}
Any complex balanced mass action system is absolutely complex balanced.
\end{theorem}
The extension is partial because while valid for arbitrary deficiency, the set of kinetics is restricted to mass action kinetics. 

This partial extension  view leads to the following question (which we call the Horn-Jackson ACB Extension Problem or simply the Extension Problem): beyond mass action kinetics, which necessary or sufficient conditions ensure absolute complex balancing in system with positive deficiency?  In particular, are there any kinetics sets other than mass action, where any complex balanced system with positive deficiency is absolutely complex balanced?

A natural candidate for extending the Horn-Jackson result is the complex balanced subset of power law kinetic systems with reactant-determined kinetic orders (denoted by PL-RDK), i.e., those where branching reactions of a reactant have identical rows in the system's kinetic order matrix. PL-RDK systems are precisely the CF power law systems and correspond to a subset of the generalized mass action systems (GMAS) introduced by S. M\"{u}ller and G. Regensburger in 2014 \cite{MURE2014}. Numerous results on mass action systems have been extended to PL-RDK systems, including the Core Deficiency Zero Theorem \cite{MURE2012,MURE2014}, toric steady states \cite{JOHN2014}, the Johnston-Siegel Linear Conjugacy Criterion \cite{CONM2018}, Birch's Theorem \cite{CMPY2019} and the Shinar-Feinberg ACR Theorem \cite{FOME2020}. However, our first main result is an example of a complex balanced PL-RDK system, which is not absolutely complex balanced, showing that the analogue of the Horn-Jackson ACB Theorem does not hold for PL-RDK systems in general. 

\begin{RE2}
Consider the weakly reversible kinetic system with three species $A_1$, $A_2$ and $A_3$ below.  
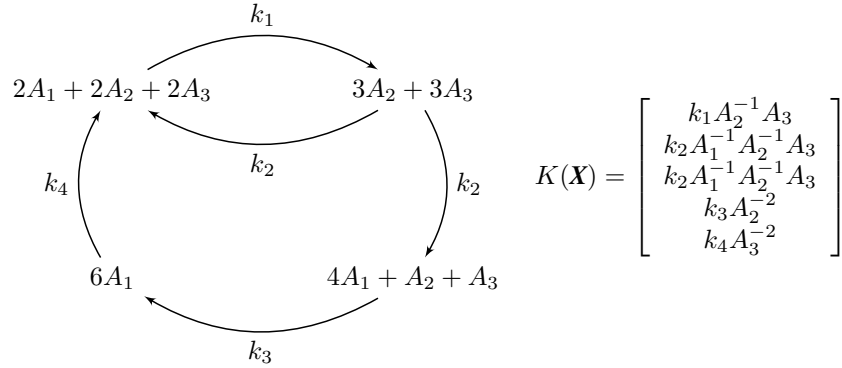
\begin{equation}
\nonumber
\begin{tikzpicture}[baseline=(current  bounding  box.center)]
\tikzset{vertex/.style = {draw=none,fill=none}}
\tikzset{edge/.style = {bend left,->,> = latex', line width=0.20mm}}
\node[vertex] (1) at  (0,0) {$2A_1 + 2A_2 +2A_3$};
\node[vertex] (2) at  (4,0) {$3A_2 + 3A_3$};
\node[vertex] (3) at  (0,-2.5) {$6A_1$};
\node[vertex] (4) at  (4,-2.5) {$4A_1+ A_2 + A_3$};
\draw [edge]  (1) to["$k_1$"] (2);
\draw [edge]  (2) to["$k_2$"] (1);
\draw [edge]  (2) to["$k_2$"] (4);
\draw [edge]  (4) to["$k_3$"] (3);
\draw [edge]  (3) to["$k_4$"] (1);
\end{tikzpicture}
\quad
K(\textbf{X})=\left[ 
\begin{array}{c}
	k_1 A_2^{-1} A_3 \\
	k_2 A_1^{-1}A_2^{-1}A_3 \\
	k_2 A_1^{-1}A_2^{-1}A_3 \\
	k_3 A_2^{-2} \\
    k_4 A_3^{-2} \\
\end{array}
 \right]
\end{equation}

Computing $\hat{T}$, we have: 

    $$\hat{T}=\left[ 
\begin{array}{cccc}
0 & -1 & 0 & 0 \\
-1 & -1 & -2 & 0 \\
1 & 1 & 0 & -2 \\
1 & 1 & 1 & 1 \\
\end{array}
 \right].$$

Hence, the system is both PL-RDK and PL-TIK since $\hat{T}$ has maximal rank. Furthermore, $\hat{\delta}=0$. By Theorem  \ref{RKDZT}, the system is complex balanced. 

We verify that the system is not ACB by a direct computation of the equilibria sets.

$$E_+(\mathscr{N},K)=\left\{ \left[ 
\begin{array}{c}
	A_1 \\
	A_2 \\
	A_3 \\
\end{array}
 \right] \in \mathbb{R}^3_+ \middle| 	k_1 A_2^{-1} A_3 - 3 k_2 A_1^{-1}A_2^{-1}A_3 - k_3 A_2^{-2} + 2k_4 A_3^{-2} = 0 \right\}.$$ 
    $$Z_+(\mathscr{N},K)=\left\{ \left[ 
\begin{array}{c}
	A_1 \\
	A_2 \\
	A_3 \\
\end{array}
 \right] \in \mathbb{R}^3_+ \middle| \begin{array}{c}
	 - k_1 A_2^{-1} A_3 + k_2 A_1^{-1}A_2^{-1}A_3 + k_4 A_3^{-2} = 0 \\
	 k_1 A_2^{-1} A_3 - 2k_2 A_1^{-1}A_2^{-1}A_3  = 0 \\
	 k_2 A_1^{-1}A_2^{-1}A_3 - k_3 A_2^{-2} = 0 \\
	k_3 A_2^{-2} - k_4 A_3^{-2} = 0 \\
\end{array} \right\}.$$ 

Consider $v \subseteq  \mathbb{R}^3$ where $$v = \left\{ \left[ 
\begin{array}{c}
	A_1 \\
	A_2 \\
	A_3 \\
\end{array}
 \right] \in \mathbb{R}^3_+ \middle| \begin{array}{c}
	\dot{A_1}  =k_1 A_2^{-1} A_3 - k_2 A_1^{-1}A_2^{-1}A_3 = 0 \\
	\dot{A_2}  = k_2 A_1^{-1}A_2^{-1}A_3 - k_3 A_2^{-2} = 0 \\
	\dot{A_3}  =-3k_2 A_1^{-1}A_2^{-1}A_3  + 2k_4 A_3^{-2} = 0 \\
\end{array} \right\}.$$ It is easy to see that the solutions of the equations in $v$ is a subset of the solutions of the equation in $E_+$ (i.e., adding the 3 equations in $v$ will result to the equation in $E_+$). Furthermore, $v \cap Z_+ = \emptyset$ since the equations in $v$ and $Z_+$ will not yield positive solutions. Now, to show $v$ is non-empty, let $k_4 = \frac{3}{2}r_4$. The set $v$ can be written as 

$$v = \left\{ \left[ 
\begin{array}{c}
	A_1 \\
	A_2 \\
	A_3 \\
\end{array}
 \right] \in \mathbb{R}^3_+ \middle| 
 \left[
 \begin{array}{cccc}
1 & -1 & 0 & 0 \\
0 & 1 & -1 & 0 \\
0 & -3 & 0 & 3 \\
\end{array} \right]
\left[ 
\begin{array}{c}
	k_1 A_2^{-1} A_3 \\
	k_2 A_1^{-1}A_2^{-1} A_3 \\
	k_3 A_2^{-2} \\
    r_4 A_3^{-2} \\
\end{array}
 \right]=0
\right\}  $$

The CRN of the ODE above can be: 
 \begin{equation}
\nonumber
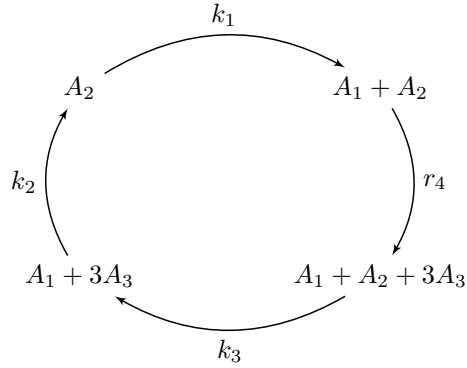
\begin{tikzpicture}[baseline=(current  bounding  box.center)]
\tikzset{vertex/.style = {draw=none,fill=none}}
\tikzset{edge/.style = {bend left,->,> = latex', line width=0.20mm}}
\node[vertex] (1) at  (0,0) {$A_2$};
\node[vertex] (2) at  (4,0) {$A_1 + A_2$};
\node[vertex] (3) at  (0,-2.5) {$A_1+3A_3$};
\node[vertex] (4) at  (4,-2.5) {$A_1+A_2+3A_3$};
\draw [edge]  (1) to["$k_1$"] (2);
\draw [edge]  (3) to["$k_2$"] (1);
\draw [edge]  (2) to["$r_4$"] (4);
\draw [edge]  (4) to["$k_3$"] (3);
\end{tikzpicture}
\end{equation}
 
 Since the network is weakly reversibly and PL-TIK, $v$ is non-empty (Theorem \ref{RKDZT}). Hence, the given kinetic system is complex balanced but not ACB.
\end{RE2}

 \begin{remark}
 Although we provide the simpler, direct proof of the example not being ACB, we were actually led to it by considerations in connection with a partial converse to the Feinberg ACB Theorem. This approach is presented in Section \ref{sec:partialconverse} as it provides a general technique to finding necessary conditions for ACB systems with positive deficiency as well as generating counterexamples.
 \end{remark}

\section{A necessary and sufficient condition for absolute complex balancing of kinetic systems of a CLP system}\label{sec:ACBCLP}

In his 1979 Wisconsin Lecture Notes, M. Feinberg showed that ACB in a complex balanced mass action system is equivalent to its set of complex balanced equilibria being ``log parametrized'' by $S^\perp$, i.e., $Z_+(\mathscr{N}, K)  = \{ x \in \mathbb R^\mathscr S_> |\log x - \log x^*| \in S^\perp\}$, where $S$ is the network's stoichiometric subspace and $x^*$ a given complex balanced equilibrium. We call a complex balanced system with this ``log parametrization'' property a CLP system with flux space $S$ (and parameter space $S^\perp$).

Thirty five years later, S. M\"{u}ller and G. Regensburger extended the theory of toric mass action systems of Craciun et al. \cite{CDSS2009} to show that any complex balanced PL-RDK system of CLP type with flux space = $\tilde{S}$ (and parameter space $\tilde{S}^\perp$), where the kinetic order subspace $\tilde{S}$, the kinetic analogue of the stoichiometric subspace $S$, is generated by the differences of kinetic complexes, which are the columns of the system's $T$ matrix, of the product and reactant complexes of the system's reactions. Hence, the example in Section \ref{sec:ACBPLRDK} also shows that the equivalence between ACB and CLP does not extend to complex balanced PL-RDK systems.

In this Section, we broaden the scope of our study to the set of CLP systems in general, which include subsets of poly-PL and Hill-type systems. We provide a necessary and sufficient condition for ACB in CLP systems, thus resolving the Extension Problem for them. In view of the previously mentioned result of M\"{u}ller and Regensburger, this also resolves the Extension Problem for complex balanced PL-RDK systems. Building on results of M. Feinberg in his 1979 Lectures, basic concepts and fundamental properties of LP (both CLP and PLP) systems are discussed in Section \ref{subsec:fundLP}.  In Section \ref{subsec:necsufCLP}, the necessary and sufficient condition is derived and illustrated with various examples of power law systems.

\subsection{Fundamentals of LP sets and LP systems}\label{subsec:fundLP}
We begin with the concept of an LP set of a chemical kinetic system $(\mathscr{N}, K)$ with reaction network $\mathscr{N}=(\mathscr{S}, \mathscr{C}, \mathscr{R})$:
\begin{definition}
An LP set is a non-empty subset of $\mathbb{R}^{\mathscr{S}}_>$ of the form $E(P,x^*):=$ where $P$ is a subspace of $\mathbb{R}^{\mathscr{S}}$ (called an LP set's flux subspace) and $x^*$ a given element of $\mathbb{R}^{\mathscr{S}}_>$ (called the LP set's reference point). $P^{\perp}$ is called an LP set's parameter subspace and the positive cosets of $P$ are called LP set's flux classes.
\end{definition}

In \cite{FEIN1979}, M. Feinberg derived the following important property of an LP set, based on the work by F. Horn and R. Jackson:

\begin{proposition}\label{LPset}
For any LP set $E=E(P,x^*)$ and any of its flux class $Q$, $|E\cap Q|=1$.
\end{proposition}
We reproduce the proof of Feinberg in a slightly more general form for the convenience of the reader. We need the following lemma:

\begin{lemma}
Let $\mathscr{S}$ be any finite set and let $\mathbb{R}^{\mathscr{S}}_>$ be the vector space generated by $\mathscr{S}$. Let $S$ be a linear subspace of $\mathbb{R}^{\mathscr{S}}_>$ and let $a$ and $b$ be elements of $\mathbb{R}^{\mathscr{S}}_>$. There exists a (unique) vector $\mu\in S^{\perp}$ such that
$$ae^\mu-b$$
is an element of $S$.
\end{lemma}

\textit{Proof.}
We first show that $|E\cap Q| >0$. Let $Q=p+S$. According to the preceding lemma, there is a unique vector $\mu\in S^{\perp}$ such that $x^*e^\mu-p\in S$. Set $x=x^*e^\mu$. Clearly, $x\in p+S=Q$. On the other hand, $\ln x=\ln x^*+\mu$ so that $\ln x\ln x^*\in S^{\perp}$. Hence $x\in E\cap Q$. We now show that $|E\cap Q|=1$.

Suppose $x,x'$ are in $E\cap Q$. Then
\begin{equation}\label{lemmaeqn1}
  0=(x'-x)\cdot(\ln x'-\ln x)=\sum_{S\in\mathscr{S}}(x'_S-x_S)(\ln x'_S-\ln x_S).  
\end{equation}
Since the function $\ln:S\rightarrow R$ is strictly monotonically increasing, \eqref{lemmaeqn1} can hold only if, for all $S\in\mathscr{S}$, $x'_S=x_S$, that is, only if $x'=x$.
$\blacksquare$

We now introduce the concepts relating LP sets and equilibria sets of kinetic systems.

\begin{definition}
A subset $E$ of the positive equilibria set $E_+(\mathscr{N}, K)$ of a chemical kinetic system $(\mathscr{N}, K)$ is of LP type if $E$ is an LP set, i.e., $E=E(P,x^*)$ for a subspace $P$ of $\mathbb{R}^{\mathscr{S}}$ and an element $x^*$ in $E$. A chemical kinetic system is of PLP type if $E_+(\mathscr{N}, K)\neq\emptyset$ and of LP type for a subspace $P_Z$ of $\mathbb{R}^{\mathscr{S}}$. An LP system is a PLP or CLP system. $P_E$ and $P_Z$ are called LP system's flux subspaces and the positive cosets in $\mathbb{R}^{\mathscr{S}}$ are the LP system's flux classes.   
\end{definition}

\begin{remark}
We note that the definition of an LP system does not depend on the given equilibrium $x^*$. If $x^{**}$ is another equilibrium, then $\log x - \log x^{**}  = \log x - \log x^* + (\log x^* - \log x^{**})$ is contained in the parameter space for any other equilibrium $x$.
\end{remark}

The following proposition justified the term "parameter subspace" for $P^{\perp}$:

\begin{proposition}\label{Pperp}
Let $(\mathscr{N}, K)$ be a chemical kinetic system.
\begin{enumerate}
    \item [(i)] If $(\mathscr{N}, K)$ is of type PLP with flux subspace $P_E$ and reference point $x^*$ in $E_+(\mathscr{N}, K)$, then the map $L_{x^*}:E_+(\mathscr{N}, K)\rightarrow P_E^\perp$ given by $L_{x^*} (x)=\log (x)-\log (x^*)$ is a bijection.
    \item [(ii)] If $(\mathscr{N}, K)$ is of type CLP with flux subspace $P_Z$ and reference point $x^*$ in $Z_+(\mathscr{N}, K)$, then the restriction to  $Z_+(\mathscr{N}, K)$ of $L_{x^*}:Z_+(\mathscr{N}, K)\rightarrow P_Z^\perp$ is a bijection.
\end{enumerate}
\end{proposition}

\textit{Proof.}
For (i), observe that $L_{x^*} (x_1)=L_{x^*} (x_2)$ if and only if $\log (x_1)=\log (x_2)$ if and only if $x_1=x_2$, due to the injectivity of $\log:\mathbb{R}^{\mathscr{S}}_>\rightarrow \mathbb{R}^{\mathscr{S}}$. For $p\in P_E^\perp$, due to the surjectivity of $\log$, there is an $x\in \mathbb{R}^{\mathscr{S}}_>$ with $\log (x)=p+\log (x^*)$. Hence $\log (x)-\log (x^*)\in P_E^\perp$. Since the kinetic system is PLP, then $x\in E_+(\mathscr{N}, K)$. The proof for (ii) is analogous to that of (i).
$\blacksquare$

\begin{example}
 Statements ii) and iii) of M. Feinberg's formulation of the Horn-Jackson ACB Theorem show that any complex balanced mass action system is a CLP system with $P_Z = S.$
\end{example}

\begin{example}
Any weakly reversible mass action systems satisfying the conditions of the Deficiency One Theorem is of PLP type.
\end{example}

\begin{example}
S. M\"{u}ller and G. Regebsburger in 2014 extended the work of Craciun et al in 2007 on toric mass action systems to establish that any complex balanced generalized mass action system (GMAS) is of CLP type. PL-RDK systems map bijectively into the set of GMAS: the only additional property  for them is that if the zero complex is a reactant, then it is mapped to the zero kinetic complex. The property ensures consistency with biochemical usage and the extension of the GMAS kinetic map (of sets) to a linear map between their linear spans. For complex balanced PL-RDK systems, $P_Z = \tilde{S}$, which is called the  kinetic order subspace.
\end{example}

\subsection{A necessary and sufficient condition for absolute complex balancing in a CLP system}\label{subsec:necsufCLP}

The next two Propositions motivate our definition of a bi-LP system.

\begin{proposition}
Let $(\mathscr{N}, K)$ be a chemical kinetic system, which is both CLP and PLP with flux subspaces $P_Z$ and $P_E$ respectively. Then $P_Z \subset P_E.$ If the system is absolutely complex balanced, then $P_Z = P_E.$
\end{proposition}

\textit{Proof.}
Let $x^*\in Z_+(\mathscr{N}, K)$. Since $Z_+(\mathscr{N}, K)$ is contained in $E_+(\mathscr{N}, K)$, the map $L_x^*$ on $Z_+(\mathscr{N}, K)$ (from Proposition \ref{Pperp}) is just the restriction of that on  $E_+(\mathscr{N}, K)$. Hence its image $P_Z$ is contained in that of  $E_+(\mathscr{N}, K)$, that is, in $P_E.$ Clearly, if the definition domains coincide, i.e. ACB holds, then $P_Z = P_E.$

Furthermore, the following statements derive directly from applying Proposition \ref{LPset} to $E = E_+(\mathscr{N}, K)$ with $P = P_E$ and $E = Z_+(\mathscr{N}, K)$ with $P = P_Z.$
$\blacksquare$

\begin{proposition}
Let $(\mathscr{N}, K)$ be a chemical kinetic system.
\begin{enumerate}
    \item [(i)] If $(\mathscr{N}, K)$ is a PLP system then $|E_+(\mathscr{N}, K)\cap Q|=1$ for any of its flux classes $Q$.
    \item [(ii)] If $(\mathscr{N}, K)$ is a CLP system then $|Z_+(\mathscr{N}, K)\cap Q|=1$ for any of its flux classes $Q$.
\end{enumerate}
\end{proposition}

Note that if $(\mathscr{N}, K)$ is both CLP and PLP with $P_Z = P_E$, then the coset intersection counts can be compared.

\begin{definition}
A kinetic system $(\mathscr{N}, K)$ is bi-LP iff it is both CLP and PLP and their flux spaces coincide, i.e. $P_Z = P_E$.
\end{definition}

The following theorem is the criterion for absolute complex balancing of a CLP system:

\begin{theorem}
\label{4.10}
Let $(\mathscr{N}, K)$ be a CLP system with flux space $P_Z$ and reference point $x^*$. Then $(\mathscr{N}, K)$ is absolutely complex balanced if and only if $(\mathscr{N}, K)$ is a bi-LP system.
\end{theorem}

\textit{Proof.}
Observe that the implication is straightforward since $(\mathscr{N}, K)$ is ACB implies it is PLP with $P_Z = P_E$ and $x^*$, hence bi-LP.

To show the converse, note that since the cosets of $\mathbb R^{\mathscr {S}} (P=P_E=P_z)$ partition $\mathbb R^{\mathscr {S}}$, the positive cosets partition $\mathbb R^{\mathscr {S}}_>$. Consequently, $E_+(\mathscr{N}, K)$ and $Z_+(\mathscr{N}, K)$ are disjoint unions of their intersections with the positive cosets. Since $E_+(\mathscr{N}, K)$ contains $Z_+(\mathscr{N}, K)$,   $E_+(\mathscr{N}, K)\cap Q$ contains $Z_+(\mathscr{N}, K)\cap Q$, the equal count assumption implies that the equilibria sets coincide.  
$\blacksquare$

\begin{example}
Statements ii)  to  v) of Feinberg's formulation implies that every complex balanced mass action system is an absolutely complex balanced CLP system with $P_Z = S$. It follows from the criterion that such a system has a unique complex balanced equilibrium in each stoichiometric class.
\end{example}

\begin{example}
Talabis et al. (2019) \cite{TAMJ2019} showed that any weakly reversible PL-TIK is unconditionally complex balanced, hence, being PL-RDK, of CLP type with $P_Z = \tilde{S}$. A weakly reversible PL-TIK system satisfying the conditions of the Deficiency One Theorem for PL-TIK in addition is of PLP type with the same flux subspace, therefore bi-LP.
\end{example}


\section{Decompositions and absolutely complex balanced systems} \label{sec:decompoACB}

In this Section and Section \ref{sec:ACBHill}, we take the first steps in addressing the Extension Problem in non-CLP systems. We develop several methods of constructing ACB systems with positive deficiency and apply them to various kinetic sets to obtain \textbf{\underline{sufficient conditions}} for their occurrence for the various kinetic types. We begin with the method of combining incidence independent and independent decompositions of the underlying network, whose subnetworks share some interesting properties.

\subsection{Decompositions and positive equilibria}

We review some concepts and results in decomposition theory and  refer to \cite{FML2021} for more details.

\begin{definition}
Let $\mathscr{N} =(\mathscr{S}, \mathscr{C}, \mathscr{R})$ be a CRN. A \textbf{covering} of $\mathscr{N}$ is a collection of subsets $\{ \mathscr{R}_1, \mathscr{R}_2,\dots, \mathscr{R}_p \}$ whose union is $\mathscr{R}$. A covering is called a \textbf{decomposition} of $\mathscr{N}$ if the sets $\mathscr{R}_i$ form a partition of $\mathscr{R}$.
\end{definition}

The definition implies that each $\mathscr{R}_i$ defines a subnetwork $\mathscr{N}_i$ of $\mathscr{N}$, namely $\mathscr{C}_i$ consisting of all complexes occurring in $\mathscr{R}_i$ and $\mathscr{S}_i$ consisting of all the species occurring in $\mathscr{C}_i$. The following proposition shows the relationship of stoichiometric subspaces induced by a covering.

\begin{proposition}[Prop. 3., \cite{FML2021}]
If $\{ \mathscr{R}_1, \mathscr{R}_2,\dots, \mathscr{R}_p \}$ is a network covering, then
\begin{enumerate}
\item[i.] $S= S_1 + S_2 + \cdots + S_p$;
\item[ii.] $s \leq s_1 + s_2 + \cdots + s_p$, 
\end{enumerate}
where $s=\dim S$ and $s_i=\dim S_i$ for $i\in \overline{1,p}$.
\end{proposition}

Feinberg  \cite{FEIN1987} defined an important class of decompositions called  independent decompositions:

\begin{definition}
A decomposition is \textbf{independent} if $S$ is the direct sum of the subnetworks' stoichiometric subspaces $S_i$ or equivalently, $s = s_1 + s_2 + \cdots + s_p$.
\end{definition}

Fortun et al. \cite{FMRL2019} derived a basic property of independent decompositions:

\begin{proposition}[Lemma 1, \cite{FMRL2019}]\label{prop:fortun}
If $\mathscr{N}=\mathscr{N}_1 \cup \mathscr{N}_2 \cup \cdots  \cup \mathscr{N}_p$ is an independent decomposition, then $
\delta \leq\delta_1 +\delta_2 + \cdots +\delta_p$, where $\delta_i$ represents the deficiency of the subnetwork $\mathscr{N}_i$.
\end{proposition}

In \cite{FEIN1987}, Feinberg found the following relationship between the positive equilibria of the ``parent network'' and those of the subnetworks of an independent decomposition.

\begin{theorem}[Rem. 5.4, \cite{FEIN1987}] \label{feinberg theorem}
Let $(\mathscr{N},K)$ be a chemical kinetic system with partition $\{\mathscr{R}_1,\mathscr{R}_2,\dots, \mathscr{R}_p \}$. If $\mathscr{N}=\mathscr{N}_1 \cup \mathscr{N}_2 \cup \cdots \cup\mathscr{N}_p$ is the network decomposition generated by the partition  and $E_+(\mathscr{N}_i,K_i)= \{ x \in \mathbb{R}^\mathscr{S}_{>0} | N_i K_i(x) = 0 \}$, then 
\begin{enumerate}
\item[i.] $ \displaystyle{\bigcap_{i\in \overline{1,p}}} E_+ (\mathscr{N}_i, K_i)\subseteq E_+ (\mathscr{N}, K)$
\item[ii.] If the network decomposition is independent, then equality holds.
\end{enumerate}
\end{theorem}

Farinas et al. \cite{FML2021} observed  the following result which led them to introduce the concept of an \textit{incidence independent decomposition}. 

\begin{proposition}[Prop. 6, \cite{FML2021}]
If $\{ \mathscr{R}_i \}$ is a network covering, then
\begin{enumerate}
\item[i.] $\text{\emph{Im} } I_a = \text{\emph{Im }} I_{a,1} + \text{\emph{Im }} I_{a,2} + \cdots + \text{\emph{Im }} I_{a,p}$, where $I_{a,i}$ denotes the incidence map of the subnetwork $\mathscr{N}_i$.
\item[ii.] $n- \ell \leq (n_1 - \ell_1) + (n_2 - \ell_2) + \cdots + (n_p - \ell_p)$, where $n-\ell=\dim I_a$ and $n_i - \ell_i=\dim I_{a,i}$ for $i\in \overline{1,p}$.
\end{enumerate}
\end{proposition}

An incidence independent decomposition was defined in \cite{FML2021} as follows: 

\begin{definition}
A decomposition $\{ \mathscr{N}_1, \mathscr{N}_2, \dots, \mathscr{N}_p \}$ of a CRN is \textbf{incidence independent} if and only if the image of the incidence map $I_a$ of $\mathscr{N}$ is the direct sum of the images of the incidence maps of the subnetworks.
\end{definition}

It follows from this definition that $n-\ell = \sum (n_i - \ell_i)$.  The linkage classes form the primary example of an incidence independent decomposition, since $n = \sum n_i$ and $\ell = \sum \ell_i$. 

The following result is the analogue of Proposition \ref{prop:fortun} for incidence independent decomposition.

\begin{proposition}[Prop. 7, \cite{FML2021}]
\label{prop:incidenceindep}
Let $\mathscr{N}=\mathscr{N}_1 \cup \mathscr{N}_2 \cup \cdots \cup \mathscr{N}_p$ be an incidence independent decomposition. Then $\delta \geq \delta_1 +\delta_2 + \cdots + \delta_p$.
\end{proposition}

A decomposition is \textbf{bi-independent} if it is both independent and incidence independent. Independent linkage class decomposition is the best known example of bi-independent decomposition.

\begin{proposition}[Prop. 9, \cite{FML2021}]
A decomposition $\mathscr{N}= \mathscr{N}_1 \cup \mathscr{N}_2 \cup \cdots \cup \mathscr{N}_p$ is independent or incidence independent and $\displaystyle{\sum_{i =1}^p} \delta_i = \delta$ if and only if $\mathscr{N}= \mathscr{N}_1 \cup \mathscr{N}_2 \cup \cdots \cup \mathscr{N}_p$ is bi-independent.
\end{proposition}

\subsection{Complex balanced systems with decompositions into ACB subnetworks}

The main result for constructing ACB systems with positive deficiency via decompositions is the following:

\begin{theorem}\label{decompoACBthm}
Let $(\mathscr{N}, K)$ be a complex balanced system. If $\mathscr{N} = \mathscr{N}_1 \cup \cdots \cup \mathscr{N}_k$ is a bi-independent decomposition into ACB subnetworks, then $(\mathscr{N}, K)$ is also absolutely complex balanced.
\end{theorem}

\textit{Proof.}
For a bi-independent decomposition, it follows from the decomposition theorems of the previous section that $E_+(\mathscr N, K) = \cap E_+(\mathscr N_i, K_i)$ and $Z_+(\mathscr N, K) = \cap Z_+(\mathscr N_i, K_i)$. Note that since $(\mathscr N, K)$ is complex balanced, all these sets are non-empty.  Since the corresponding sets in the intersections are equal, the claim follows. Note also that, since in a bi-independent decomposition, $\delta = \delta_1+\dots+ \delta_k$, a positive subnetwork deficiency implies a positive network deficiency, too.
$\blacksquare$

The previous theorem leads in a special case to a corollary of the following recent result of L. Fontanil and E. Mendoza \cite{FOME2021}: 

\begin{theorem}\label{Fontanil}
Let $(\mathscr{N}, K)$ be a weaky reversible power law system with a complex balanced PL-RDK decomposition $\mathscr{D}: \mathscr{N} = \mathscr{N}_1 \cup \cdots \cup \mathscr{N}_k$ with $P_{Z,i}=\tilde{S_i}$. If $\mathscr{D}$ is incidence independent and the induced covering $\tilde{\mathscr{D}}$ is independent, then $(\mathscr{N}, K)$ is a weakly reversible CLP system with $P_Z=\sum \tilde{S_i}$.
\end{theorem}

\begin{corollary}
Let $(\mathscr{N}, K)$ satisfy the assumptions of Theorem \ref{Fontanil}. In addition, let the decomposition be independent and its PL-RDK subnetworks be ACB. Then $(\mathscr{N}, K)$ is bi-LP with flux space sum $\sum \tilde{S_i}$.
\end{corollary}

\textit{Proof.}
By Theorem \ref{Fontanil}, $(\mathscr{N}, K)$ is CLP with flux space $\sum \tilde{S_i}$i, hence it is necessarily complex balanced. Applying Theorem \ref{decompoACBthm}, it follows that it is ACB. Theorem \ref{4.10} then implies the claim.
$\blacksquare$

\begin{remark}
The Corollary also follows from combining Theorem \ref{Fontanil}. with a recent result of B. Hernandez and E. Mendoza \cite{HEME2021} after applying Theorem \ref{4.10} to the subnetworks with flux spaces $\tilde{S_i}$.
\end{remark}

For the general case, i.e., arbitrary kinetics, in the previous theorem, we need the strong hypothesis of bi-independence for the decomposition. We next present a special construction for a class of power law systems derived from poly-PL kinetics, which needs only the weaker property of incidence independence.

\subsection{A brief review of poly-PL systems and the STAR-MSC transformation}

Poly-PL kinetic systems (denoted as PYK systems) were introduced by Talabis et al. \cite{TMMJ2020} and applied to reaction network models of evolutionary games as proposed by \cite{VELO2014}. PL-representations were used to show that a subset of weakly reversible CF poly-PL systems denoted as PY-TIK systems possessed unconditional complex balancing (UCB), i.e., the existence of a complex balanced equilibrium for any set of rate constants. The STAR-MSC transformation from a poly-PL system to a dynamically equivalent PL system was first studied in Magpantay et al \cite{MHRM2020} and used to provide a computational approach to determining whether a poly-PL system had the capacity for multistationarity, i.e. the occurrence of two or more positive equilibria in a stoichiometric class for some rate constants.

\begin{definition}
A poly-PL kinetics is a kinetics $K: \Omega \rightarrow \mathbb{R}^{\mathscr{R}}$ such that for each reaction $q$

$$K_q(x) = k_q(a_{q,1} x^{F_{q,1}} + \cdots  + a_{q,h_q}x^{ F_{q,h_q}})$$

with $k_j > 0$, $a_{q,j} > 0$, $F_{q,j} \in \mathbb{R}^{\mathscr{S}}$ and $j = 1,\cdots ,h_q$. If $h = \max h_j$, we normalize the length of each kinetics to $h$ by replacing $a_{q,h_q}x^{F_{q,h_q}}$ with $(h - h_j + 1)$ terms each equal to $\frac{1}{(h - h_j +1)}a_{q,h_q}$. We call $h$ the length of the poly-PL kinetics. The set of poly-PL kinetics of length $h$ is denoted by $PYK_h$ and $PYK = \cup PYK_h$.
\end{definition}

By subsuming the coefficients $a_{q,j}$ into new rate constants $k_{q,j} := k_q a_{q,j}$, we can write $K_q(x) = K_{q,1} +\cdots + K_{q,h}$. This relationship defines a PL-representation of the poly-PL kinetics $K = K_1 +\cdots + K_h$, with $K_j = k_j x^{F_j}$ being PLK and the rows of $F_j$ are the $F_{q,j}$. It is sometimes convenient to rearrange the kinetics into lexicographic order before normalizing the length. In that case, we call the resulting representation the canonical PL- representation of the poly-PL kinetics. 

\begin{RE3}
\sloppy
Consider a reversible Michaelis-Menten kinetics from Enzyme biology and the network with set of biochemical species $S = \left\{S_1, S_2, S_3, S_4\right\}$, set of complexes $\left\{S_1+S_2, S_1+S_3, S_4\right\}$, and reactions:
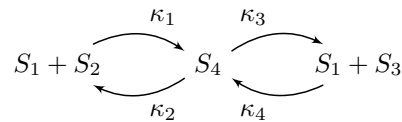
\begin{equation}
\nonumber
\begin{tikzpicture}[baseline=(current  bounding  box.center)]
\tikzset{vertex/.style = {draw=none,fill=none}}
\tikzset{edge/.style = {bend left,->,> = latex', line width=0.20mm}}
\node[vertex] (1) at  (0,0) {$S_1+S_2$};
\node[vertex] (2) at  (2,0) {$S_4$};
\node[vertex] (3) at  (4,0) {$S_1+S_3$};
\draw [edge]  (1) to["$\kappa_1$"] (2);
\draw [edge]  (2) to["$\kappa_2$"] (1);
\draw [edge]  (2) to["$\kappa_3$"] (3);
\draw [edge]  (3) to["$\kappa_4$"] (2);
\end{tikzpicture}
\end{equation}
This network assumes the Michaelis-Menten enzyme mechanism, in which a substrate $S_2$ is modified into a substrate $S_3$ through the formation of an intermediate $S_4$. The reaction is catalyzed by an enzyme $S_1$. To study $Z_+(\mathscr{N},K)$ and $E_+(\mathscr{N},K)$, Fortun et al. \cite{FTJM2020} consider the poly-PL kinetics $K(x)$ with the canonical form:
\begin{equation}
\nonumber
K(x)=\left[ 
\begin{array}{c}
k_1 (S_1 S_2+k_3 S_1 S_2^2 + k_4 S_1 S_2 S_4) \\
k_2 (S_3+k_3 S_2 S_3 + k_4 S_3 S_4) \\
k_1  (\frac{1}{3}S_2 + \frac{1}{3}S_2 + \frac{1}{3}S_2) \\
k_2  (\frac{1}{3}S_4 + \frac{1}{3}S_4 + \frac{1}{3}S_4) \\
\end{array}
 \right].
\end{equation}
$K_1(x)$, $K_2(x)$ and $K_3(x)$ are
\begin{equation}
\nonumber
K_1(x)=\left[ 
\begin{array}{c}
k_1 S_1 S_2 \\
k_2 S_3 \\
 \frac{k_1}{3} S_2 \\
\frac{k_2}{3}S_4 \\
\end{array}
 \right] \quad K_2(x)=\left[ 
\begin{array}{c}
k_1  k_3 S_1 S_2^2 \\
k_2 k_3 S_2 S_3 \\
\frac{k_1}{3} S_2  \\
\frac{k_2}{3}S_4 \\
\end{array}
 \right] \quad K_3(x)=\left[ 
\begin{array}{c}
 k_1  k_4 S_1 S_2 S_4 \\
 k_2 k_4 S_3 S_4 \\
 \frac{k_1}{3}S_2 \\
 \frac{k_2}{3}S_4 \\
\end{array}
 \right],
\end{equation}
respectively.
\end{RE3}
 The following proposition is easy to verify:
\begin{proposition}
A poly-PL kinetics $K$ is complex factorizable $\Leftrightarrow$ for any two reactions $q$, $q'$ with the same reactant complex, its canonical PL-representation $K = K_1 +\cdots + K_j$, the $(\mathscr{N},K_j)$ are PL-RDK and for its positive coefficients $a_{q,j} = a_{q',j}$, for $j = 1,\cdots ,k$.
\end{proposition}

\textit{Proof.}
The interaction functions $I_{K,q}$ and $I_{K,q'}$ are equal $\Leftrightarrow$ the corresponding power law exponents and positive coefficients coincide.
$\blacksquare$

\begin{corollary}
A poly-PL kinetics $K$ is non-complex factorizable (NF) iff for at least one $j$, $(\mathscr{N}, K_j)$ is PL-NDK.
\end{corollary}

We introduce the PYK subset needed for the construction:

\begin{definition}
 A poly-PL system $(\mathscr{N},K)$ is \textbf{PL-equilibrated (PL-complex balanced)} if  $E_+(\mathscr{N},K) = E_{+,PLE}(\mathscr{N},K)$ $\left(Z_+(\mathscr{N},K) = Z_{+,PLC}(\mathscr{N},K) \right)$ where
$$
E_{+,\text{PLE}} (\mathscr{N},K) =\displaystyle{\bigcap_{j \in \overline{1,h}}} E_+(\mathscr{N},K_j) \left(Z_{+,\text{PLC}} (\mathscr{N},K) =\displaystyle{\bigcap_{j \in \overline{1,h}}} Z_+(\mathscr{N},K_j) \right).$$
\end{definition}

The examples of PL-equilibrated (PL-complex balanced) poly-PL systems are PL-independent (PL-incidence independent) systems (s. \cite{FTJM2020} for details). 

The S-invariant termwise addition of reactions via maximal stoichiometric coefficients (STAR-MSC) method is based on the idea to use the maximal stoichiometric coefficient (MSC) among the complexes in the CRN to construct reactions whose reactant complexes and product complexes are different from existing ones. This is done by uniform translation of the reactants and products to create a ``replica'' of $\mathscr{N}$. The method creates $h - 1$ replicas of the original network and hence its transform, $\mathscr{N}^*$ becomes the union (in the sense of \cite{GHMS2020}) of the replicas and the original CRN. 

We now describe the STAR-MSC transformation. Since the domain of definition of a poly-PL kinetics is $\mathbb{R}^{\mathscr{S}}_{>}$, all $x = (X_1,X_2,\cdots ,X_m)$ are positive vectors. Let $M = 1 + \max \left\{ y_i | y \in \mathscr{C} \right\}$, where the second summand is the maximal (positive integer) stoichiometric coefficient.

For any positive integer $z$, define the vector $z$ to be the vector $(z, z,\cdots, z) \in R^{\mathscr{S}}$. For each complex $y \in \mathscr{C}$, form the $(h-1)$ complexes

$$y+\underline{M}, y + \underline{2M},\cdots , y + \underline{(h - 1)M}.$$

Each of these complexes are different from all existing complexes and each other as shown in the following proposition:

\begin{proposition}
Let $\mathscr{N}^*  =(\mathscr{S},\mathscr{C}^*,\mathscr{R}^*)$ be the STAR-MSC transform of $\mathscr{N}=(\mathscr{S},\mathscr{C},\mathscr{R})$, $\mathscr{N}^*_1 := \mathscr{N}$ and $\mathscr{N}^*_j$ is the subnetwork defined by $\mathscr{R}^*_{j-1}$ for $j \in \overline{2, h}$. Then $|\mathscr{C}^*| = h \cdot n$ and $|\mathscr{R}^*| = h \cdot r$. 
\end{proposition}

\subsection{STAR-MSC transformation of weakly reversible deficiency zero poly-PL systems}
The last proposition in the previous section immediately implies the following result:

\begin{proposition}
The deficiency of $(\mathscr{N}^*, K^*)$ is 
$$\delta^*= \delta + (n - l)(h - 1).$$
\end{proposition}

\textit{Proof.}
$\delta^* = n^* - l^* - s^* = nh - lh - s = (n - l - s) + (n - l)(h - 1) = \delta + (n - l)(h - 1)$.
$\blacksquare$

The following result shows the claimed properties of the constructed class:

\begin{theorem}
 The STAR-MSC transform $(\mathscr{N}^*, K^*)$ of a weakly reversible, deficiency zero, PL-complex balanced PY-RDK system of length $h \geq 2$ has positive deficiency, is non-PL-FSK and absolutely complex balanced.
 \label{taas}
\end{theorem}

\textit{Proof.}
By STAR-MSC construction, $\mathscr{N}^*$ is weakly reversible.  From the previous Proposition, $\delta^*= 0 + (n - l)(h - 1) > 1(2-1) = 1$, since $n - l > 1$. The kinetic order matrix of $K^*$ consists of at least two stacked copies of $F_1$ implies that $K^*$ is not factor span surjective. As shown in \cite{FTJM2020}, $Z_+(\mathscr{N}^*,K^*) = Z_{+,PLC}(\mathscr{N},K)$, since $(\mathscr{N},K)$ is PL-complex balanced, $Z_{+,PLC}(\mathscr{N}, K) = Z_+(\mathscr{N}, K)$. $\mathscr{N}$ has a CB equilibrium by assumption and zero deficiency yielding $Z_+(\mathscr{N},K) = E_+(\mathscr{N},K)$. Therefore, $(\mathscr{N}^*,K^*)$ is dynamically equivalent to $(\mathscr{N},K)$, $E_+(\mathscr{N}, K) = E_+(\mathscr{N}^*,K^*)$, which proves the claim.
$\blacksquare$

\subsection{The decomposition $\mathscr{N}^* =\mathscr{N}^*_1 \cup \cdots \cup \mathscr{N}^*_h$}

By assumption, $\mathscr{N}$ is complex balanced and PL-complex balanced. Hence, the $\mathscr{N}_j$ are complex balanced. Since the $\mathscr{N}^*_j$ are replicas of $\mathscr{N}$ and have the same kinetics as $\mathscr{N}_j$, they are also complex balanced  and have zero deficiency, hence, they are ACB. The decomposition above is the linkage class decomposition of $\mathscr{N}^*$, hence $\mathscr{N}^*$ is complex balanced. The decomposition is incidence independent , but not bi-independent (because $\delta^* > 0 = $ sum of subnetwork deficiencies). Therefore, $\mathscr{N}^*$ is an example of a complex balanced PLK system which has a decomposition into ACB subnetworks, though only incidence independent, is nevertheless ACB. If any of the $(\mathscr{N}_j, K_j)$ is PL-NDK, then $(\mathscr{N}^*, K^*)$ is also PL-NDK.

\section{Absolutely complex balanced Hill-type and weakly monotonic kinetic systems}\label{sec:ACBHill}

In this Section, we continue our description of methods for constructing ACB systems of positive deficiency for sets of kinetics which may be non-CLP. Section 6.1 introduces the Positive Function Factor (PFF) method which in Section 6.2 is applied to the study of ACB in Hill type systems. In Section 6.3 currently known results about ACB in poly-PL systems are presented. Section 6.4 introduces the Coset Intersection Count (CIC) method and applies it in Section 6.5 to complex balanced weakly monotonic kinetic systems on conservative and concordant networks.

\subsection{The Positive Function Factor (PFF) Method}
Let $(\mathscr{S},\mathscr{C},\mathscr{R})$ be a CRN and $\mathscr{K}_\Omega (\mathscr{N})$ be the set of kinetics on $\mathscr{N}$ and definition domain $\Omega$. $\mathscr{K}_\Omega (\mathscr{N})$ is a commutative semi-ring with respect to component wise addition and multiplication. We introduce the following equivalence relation in $\mathscr{K}_\Omega (\mathscr{N})$:

\begin{definition}
Two kinetics $K, K'$ in $\mathscr{K}_\Omega (\mathscr{N})$ are positive function factor equivalent (\textbf{PPF-equivalent}) if for all $x\in \mathbb{R}^{\mathscr{S}}_{>}$ and every reaction $q$, $\dfrac{K_q(x)}{K'_q(x)}$ is a positive function of $x$ only, i.e. independent of $q.$
\end{definition}

It is easy to see that this relation is indeed an equivalence relation in $\mathscr{K}_\Omega (\mathscr{N})$. A key property of PFF-equivalence is expressed in the following Proposition:

\begin{proposition}
If $K$, $K'$ are PPF, then
\begin{enumerate}
    \item [i.] $Z_+(\mathscr{N}, K) =  Z_+(\mathscr{N}, K')$
    \item [ii.] $E_+(\mathscr{N}, K) =  E_+(\mathscr{N}, K')$
    \item [iii.] $(\mathscr{N},K)$ has ACB (i.e. absolutely complex balanced) $\Leftrightarrow$ $(\mathscr{N},K')$ has ACB. 
\end{enumerate}
\end{proposition}

\textit{Proof.}
For a reaction $q: y_q \rightarrow y'_q$, we denote the characteristic functions of the reactant and product complexes by $\omega_q$ and $\omega_q'$, respectively. To show (i), consider $K(x)$'s CFRF (complex formation rate function) and we have $g(x)=\sum_{q \in \mathscr{R}}k_qK_q(x) (\omega_q' -\omega_q) = \sum_{q \in \mathscr{R}}k_qU(x)K'_{q}(x) (\omega_q' -\omega_q) = U(x) \sum_{q \in \mathscr{R}}k_qK_q(x) (\omega_q' -\omega_q)$, where $U(x)$ is the positive function factor. Hence $g(x) = U(x)g'(x)$, and since $U(x)>0$ for $x>0$, $g(x)=0$ if and only if $g'(x)=0$ for such $x$, which shows i). A similar argument shows ii). Clearly, iii) follows from i) and ii). 
$\blacksquare$

\begin{remark}
\begin{enumerate}
    \item [1.] If the function $U(x)$ is constant, i.e. $U(x) = a$, then PFF-equivalence leads to linear conjugacy with conjugacy vector $(a,...,a)$.
    \item [2.] The types of the kinetics $K$, $K'$ can determine the type of the function $U(x)$. For example, if both $K$, $K'$ are power law kinetics, then $U(x)$ is also a power law function.
    \item [3.] In Section \ref{sec:ACBpolyHILL}, we discuss PFF-equivalence of Hill type and various enzymatic kinetics with sums of power law kinetics.
\end{enumerate}
\end{remark}

Statements i) and ii) suggest the following further relation in $K_\Omega(\mathscr{N})$:

\begin{definition}
Two kinetics $K$, $K'$ in $K_\Omega(\mathscr{N})$ are \textbf{equilibria sets coincident} (ESC) iff $Z_+(\mathscr{N}, K) =  Z_+(\mathscr{N}, K')$ and $E_+(\mathscr{N}, K) =  E_+(\mathscr{N}, K')$.
\end{definition}

As already observed above, for equilibria sets coincident $K$, $K'$, $(\mathscr{N}, K)$ is ACB $\Leftrightarrow$ $(\mathscr{N}, K')$ is ACB. In the next Section, we extend the relation to kinetic systems $(\mathscr{N}, K)$ and $(\mathscr{N}^*, K^*)$ where $\mathscr{N}$ and $\mathscr{N}^*$ have the same set of species, but different sets of complexes and reactions. The ESC relation also is not induced by positive function factor relationship.

\subsection{Absolute complex balancing in Hill type kinetic systems}
\label{sec:ACBpolyHILL}
In \cite{WIFE2013}, Wiuf and Feliu introduced the set of Hill type kinetics (HTK) as follows:

\begin{definition}
\sloppy
A \textbf{Hill-type kinetics} assigns to each $q^{th}$ reaction, with $q=1,2,\ldots,r$, a function $K_q: \mathbb{R}^{\mathscr{S}}_{\geq} \rightarrow \mathbb{R}$ of the form
\[ K_q(x)=k_q \frac{\prod_{i=1}^m x_i^{F_{qi}}}{\prod_{i=1}^m (d_{qi}+x_i^{F_{qi}})}\]
for $i=1,\ldots,m$, where the rate constant $k_q > 0$, $F:=(F_{qi})$ and $D:=(d_{qi})$ are $r \times m$ real and nonnegative real matrices called the kinetic order matrix and dissociation constant matrix, respectively. Furthermore, $\supp{D_q}=\supp{F_q}$ to ensure normalization of zero entries to 1.  
\end{definition}

Rate functions of this type studied by A. Hill in 1910 for the case of one species (i.e., $m=1$) and non-negative integer exponents \cite{H1910}. Several years later, L. Michaelis and M. Menten focused their investigations in 1913 on functions with exponent 1 \cite{MM1913}. Refinements and extensions of the Michaelis-Menten model were widely applied to enzyme kinetics \cite{S1975} and also found their way into models of complex biochemical networks in Systems Biology. \\

In \cite{HEME2021}, Hernandez and Mendoza defined for any Hill type kinetic system $(\mathscr{N}, K)$ an associated poly-PL kinetic system $(\mathscr{N}, K_{PY})$ with the following key property:

\begin{theorem}
For any HTK system $(\mathscr{N}, K)$ we have:
\begin{enumerate}
    \item[(i)] $Z_+(\mathscr{N}, K)=Z_+(\mathscr{N}, K_{PY})$
    \item[(ii)] $E_+(\mathscr{N}, K)=E_+(\mathscr{N}, K_{PY})$
    \item[(iii)] Like $K$, $K_{PY}$ is defined on the whole non-negative orthant.
\end{enumerate}
\end{theorem}

The following Corollary follows directly from the equations in statements (i) and (ii):

\begin{corollary}
$(\mathscr{N}, K)$ is absolutely complex balanced $\Longleftrightarrow$ $(\mathscr{N}, K_{PY})$ is absolutely complex balanced.
\end{corollary}

\begin{RE3}
Fortun et al. \cite{FTJM2020} also consider the associated rational type kinetics $K_{Q}(X)$ of Running example 3, it is given by
$$K_{Q}(X) = \left[ 
\begin{array}{c}
k_1 S_1 S_2\\
k_2  S_3 \\
k_1  \frac{S_2}{1+k_3 S_2 + S_4} \\
k_2  \frac{S_4}{1+k_3 S_2 + S_4} \\
\end{array}
 \right].$$
 Since $(\mathscr{N},K)$ is absolutely complex balanced, so is $(\mathscr{N},K_Q)$.
\end{RE3}

Hence, the study of absolute complex balancing in poly-PL kinetic systems will also yield results for the property in Hill type kinetic systems. In Section 6 of \cite{HEME2021},  the authors show that an analogous association of a poly-PL system also holds for a proper superset of Hill type kinetic systems, which includes many other rate functions used in enzymatic systems. Since the equations of statements (i) and (ii) also hold, then the equivalence of ACB in such systems to that in the associated poly-PL systems is also valid.

\subsection{Absolute complex balancing in PL-complex balanced PYK kinetic systems}

As an initial result on absolute complex balancing in poly-PL systems, we provide a characterization of the property in PL-complex balanced systems. We begin with a concept naturally suggested by the containment of \\ $Z_{+,PLC}(\mathscr{N},K)$ in $E_{+,PLE}(\mathscr{N},K)$ where $E_{+,PLE} (\mathscr{N},K)$ and $Z_{+,PLC} (\mathscr{N},K)$ are defined  as 
$$
E_{+,PLE} (\mathscr{N},K) =\displaystyle{\bigcap_{j \in \overline{1,h}}} E_+(\mathscr{N},K_j) \text{ and } Z_{+,PLC} (\mathscr{N},K) =\displaystyle{\bigcap_{j \in \overline{1,h}}} Z_+(\mathscr{N},K_j).
$$

\begin{definition}
\sloppy
A weakly reversible poly-PL system $(\mathscr{N},K)$ is \textbf{absolutely PL-complex balanced} iff $Z_{+,PLC}(\mathscr{N},K) \neq \emptyset$ implies $Z_{+,PLC}(\mathscr{N},K)=E_{+,PLE}(\mathscr{N},K)$.
\end{definition}

We have the following characterization of absolute complex balancing in a weakly reversible PL-complex balanced system:

\begin{theorem}
Let $(\mathscr{N},K)$ be a weakly reversible, PL-complex balanced poly-PL system. Then the following statements are equivalent:
\begin{enumerate}
    \item[(i)] $(\mathscr{N},K)$ is absolutely complex balanced.
    \item[(ii)] $(\mathscr{N}^*,K^*)$ is absolutely complex balanced.
    \item[(iii)] $(\mathscr{N},K)$ is PL equilibrated and is absolutely PL complex balanced.
\end{enumerate}
\end{theorem}

\textit{Proof.}
For $(i) \Longleftrightarrow (ii)$. Since $\mathscr{N}_j^*$ is a replica of $\mathscr{N}$ and $K_j^* = K_j$ on the corresponding reactions, we have $Z_+(\mathscr{N},K_j) = Z_+(\mathscr{N}_j^*,K_j^*)$ and hence \[ Z_{+,PLC} (\mathscr{N},K)=Z_{+}(\mathscr{N}^*,K^*).\] Since the system is PL-complex balanced, $Z_+(\mathscr{N},K) = Z_+(\mathscr{N}^*,K^*)$ and since STAR-MSC is a dynamic equivalence, $E_+(\mathscr{N},K) = E_+(\mathscr{N}^*,K^*)$. ACB of $(\mathscr{N},K)$ is equality of the LHS, ACB of $(\mathscr{N}^*,K^*)$ equality of the RHS. \\

For $(i) \Longleftrightarrow (iii)$. Note that we have the following ``rectangle'' of containments.

\begin{enumerate}
    \item[(I)] $Z_{+,PLC}(\mathscr{N},K) \subset E_{+,PLE}(\mathscr{N},K)$
    \item[(II)] $Z_{+,PLC}(\mathscr{N},K) \subset Z_{+}(\mathscr{N},K)$
    \item[(III)] $E_{+,PLE}(\mathscr{N},K) \subset E_{+}(\mathscr{N},K)$
    \item[(IV)] $Z_{+}(\mathscr{N},K) \subset E_{+}(\mathscr{N},K)$
\end{enumerate}

Equality in (I) (PL-ACB) and (III) (PL-equilibrated) imply the equality of 
\begin{equation}
    \label{7.2}
    Z_{+,PLC}(\mathscr{N},K) = E_{+}(\mathscr{N},K).
\end{equation}

Since $Z_{+}(\mathscr{N},K)$ lies between both sets, we obtain $(\mathscr{N},K)$ is PL-complex balanced and ACB. $(\Longleftarrow)$, Equality in II (PL-complex balanced) and IV (ACB) imply Equation \ref{7.2}, too. Since $E_{+,PLE}(\mathscr{N},K)$  lies in between, we obtain $(\mathscr{N},K)$ PL-equilibrated and PL-ACB.
$\blacksquare$

We also obtain a sufficient condition for absolute complex balancing in a weakly reversible PL equilibrated poly-PL system:

\begin{corollary}
Let $K$ be a weakly reversible PL equilibrated system. If all the terms in a PL representation $K = K_1 + \hdots + K_h$ are absolutely complex balanced, then $K$ is absolutely complex balanced.
\end{corollary}

\textit{Proof.}
Since $Z_+(\mathscr{N},K_j) \neq \emptyset \Rightarrow Z_+(\mathscr{N},K_j)  =  E_+(\mathscr{N},K_j)$ for each $j = 1,\hdots,h$, $Z_{+,PLC}(\mathscr{N},K) \neq \emptyset \Rightarrow Z_{+,PLC}(\mathscr{N},K_j)  =  E_{+,PLE}(\mathscr{N},K)$, or the systems is PL-ACB. Since $K$ is PL-equilibrated, the last term is equal to $E_+(\mathscr{N},K)$. This implies that the system is PL complex balanced, and according to the Theorem, absolutely complex balanced.
$\blacksquare$

\subsection{The Coset Intersectiom Count (CIC) Method}
For a subspace $W$ of $\mathbb{R}^{\mathscr{S}}$, an element $\textbf{Q}$ of the quotient space $\mathbb{R}^{\mathscr{S}}/ W$ is called a positive coset if $\textbf{Q} \cap \mathbb{R}^{\mathscr{S}}_{\geq} \neq \emptyset.$ The Coset Intersection Count (CIC) method is based on the following Proposition:

\begin{proposition}
Let $\left( \mathscr{N},K \right)$ be a kinetic system. If $0<|E_+(\mathscr{N},K)\cap \textbf{Q}|=|Z_+(\mathscr{N},K)\cap \textbf{Q}|<\infty$  for every positive coset of a subspace $W$ of $\mathbb{R}^{\mathscr{S}}$, then $\left(\mathscr{N},K \right)$ is absolutely complex balanced.
\end{proposition}

\textit{Proof.}
Since the cosets in $\mathbb{R}^{\mathscr{S}}/ W$ partition $\mathbb{R}^{\mathscr{S}}$, the positive cosets partition $\mathbb{R}^{\mathscr{S}}_{\geq}$. Consequently,  $E_+(\mathscr{N},K)$ and $Z_+(\mathscr{N},K)$ are disjoint unions of their intersections with the positive cosets. Since $E_+(\mathscr{N},K)\cap Q$  contains $Z_+(\mathscr{N},K)\cap Q$, the equal count assumption implies that the equilibria sets coincide. 
$\blacksquare$

Most ``coset intersection count'' results in the CRNT literature involve uniqueness of the element in the intersection. The best known are those for weakly reversible mass systems with low deficiency ($\delta = 0$ or 1) with respect to $W = S$ (stoichiometric subspace of the CRN). Our example below also originates from corresponding ``low deficiency'' results for a class of power law systems. 

The CIC method can be formulated as follows: given a set of kinetic systems with a known count for $|E_+(\mathscr{N},K)\cap Q|$ or $|Z_+(\mathscr{N},K)\cap Q|$ for positive cosets of a subspace $W$, the method attempts to identify a subset with the complementary intersection count, and verify equality for each Q. Inequality for at least one Q of course implies that the system is not ACB.

\subsection{Absolute complex balancing in weakly monotonic kinetic systems}

We begin by identifying an interesting subset of complex balanced systems:

\begin{definition}
A complex balanced system $(\mathscr{N},K)$ is stoichiometrically complex balanced (SCB) iff $Z_+(\mathscr{N},K)$ has a non-empty intersection with each stoichiometric class.
\end{definition}

SCB systems include: i) any open complex balanced system, ii) any complex balanced mass action system, which has a unique complex balanced equilibrium in each stoichiometric class and PL-RDK systems satisfying the conditions of the following
Theorem of Craciun et al [4]:

\begin{theorem}
Let $(\mathscr{N}, K)$ be a complex balanced PL-RDK system with $\dim S= \dim \widetilde{S}$ and for the sign spaces of $S$ and $\widetilde{S}$, $\sigma(S)\subseteq\overline{\sigma(\widetilde{S})}$. Then $|Z_+(\mathscr{N},K) \cap Q| = 1$ for each positive stoichiometric class $Q$.
\end{theorem}

Let $L:R^\mathscr{S}\rightarrow S$ be the map defined by
$$L\alpha=\sum_{y\rightarrow y'\in\mathscr{R}}\alpha_{y\rightarrow y'} (y'-y).$$

\begin{definition}
The reaction network $\mathscr N=(\mathscr{S}, \mathscr{C}, \mathscr{R})$ is \textbf{concordant} if there do not exist an $\alpha\in\ker L$ and a nonzero $\sigma\in S$ having the following properties:
\begin{itemize}
    \item [(i)] For each $y\rightarrow y'\in\mathscr{R}$ such that $\alpha_{y\rightarrow y'}\neq 0$, $\supp{y}$ contains a species $s$ for which sgn $\sigma_s=$ sgn $\alpha_{y\rightarrow y'}$.
    \item [(ii)] For each $y\rightarrow y'\in\mathscr{R}$ such that $\alpha_{y\rightarrow y'}=0$, either $\sigma_s=0$ for all $s\in\supp{y}$ or else $\supp{y}$ contains species $s$ and $s'$ for which sgn $\sigma_s=-$ sgn $\sigma_s'$ both not zero.
\end{itemize}
\end{definition}
A network that is not concordant is discordant.

The concordance of a network is closely related to the following set of kinetics on it:

\begin{definition}
A kinetics $K$ for reaction network $\mathscr N=(\mathscr{S}, \mathscr{C}, \mathscr{R})$ is \textbf{weakly monotonic} if, for each pair of compositions $c^*$ and $c^{**}$, the following implications hold for each reaction $y\rightarrow y'\in\mathscr{R}$ such that $\supp{y}\subset\supp{c^*}$ and $\supp{y}\subset\supp{c^{**}}$:
\begin{itemize}
    \item [(i)] $K_{y\rightarrow y'}(c^{**})>K_{y\rightarrow y'}(c^{*})$ implies there is a species $s\in\supp{y}$ with $c^{**}_s>c^{*}_s$;
    \item [(ii)] $K_{y\rightarrow y'}(c^{**})=K_{y\rightarrow y'}(c^{*})$ implies $c^{**}_s=c^{*}_s$ for all $s\in\supp{y}$ or else there are species $s,s'\in\supp{y}$ with $c^{**}_s>c^{*}_s$ and $c^{**}_{s'}<c^{*}_{s'}$.
\end{itemize}
\end{definition}

For example, the set of non-inhibitory kinetics PL-NIK, i.e. those whose kinetic orders are all non-negative, constitute the weakly monotonic subset of power law kinetics.

Recall that a network is said to have injectivity in a set of kinetics if every kinetics from the set is injective, i.e. the SFRF is an injective map.  We have the following characterization of concordance in terms of the set of weakly monotonic kinetics on it:

\begin{theorem}
A reaction network has injectivity in all weakly monotonic kinetic systems derived from it if and only if the network is concordant.
\end{theorem}

The following result of Shinar and Feinberg derives the existence and uniqueness of equilibria  (``the dull behavior'') for a large class of kinetics on a large class of concordant networks:

\begin{theorem}
\label{T4.15}
If $K$ is a continuous kinetics for a conservative reaction network $\mathscr N=(\mathscr{S}, \mathscr{C}, \mathscr{R})$, then the kinetic system $(\mathscr{N},K)$ has an equilibrium within each stoichiometric compatibility class. If the network is weakly reversible and concordant, then within each nontrivial stoichiometric compatibility class there is a positive equilibrium. If, in addition, the kinetics is weakly monotonic, then that positive equilibrium is the only equilibrium in the stoichiometric compatibility class containing it.  
\end{theorem}

For example, in the context of the Coset Intersection Count method, the last statement of the Theorem can be expressed as follows: for a PL-NIK element on a weakly reversible, conservative and concordant network and any positive stoichiometric class $Q$, $|E_+(\mathscr{N},K) \cap Q| = 1.$

The result above does not guarantee that, even if the system is complex balanced, that there is a complex balanced equilibrium in every positive stoichiometric class.

If in the Theorem of Shinar and Feinberg above, the weakly monotonic kinetics is stoichiometrically complex balanced,  then the CIC allows the conclusion that the system is absolutely complex balanced.

\section{A partial converse to Feinberg's Deficiency Zero ACB Theorem}\label{sec:partialconverse}

The Horn-Jackson ACB Theorem can be viewed as providing a class of counterexamples to the full converse of the Feinberg ACB Theorem. In this Section, we derive a partial converse by identifying the set of KSE kinetics. This result leads to a \textbf{\underline{necessary}} condition for the occurrence of an ACB system with positive deficiency (for any kinetics) as well as a technique for generating complex balanced but non-ACB systems (if possible) for a particular set of kinetics. We illustrate the latter with how we discovered the example in Section 3. We begin with a brief discussion of POR kinetics, which provide a convenient way of computing dimensions in the context of KSE systems.

\subsection{Positive orthant restricted (POR) kinetics}
In CRNT, a general kinetics has usually been defined on the whole nonnegative orthant (s. \cite{FEIN2019}, Def. 3.2.1). However, for power law kinetics, which have been extensively used in modeling biochemical systems (e.g. see \cite{VOIT2013}), it is essential to exclude zero values for species involved in an inhibitory interaction, i.e. have negative kinetic orders. Accordingly, Wiuf and Feliu \cite{WIFE2013} introduced kinetics defined on subsets $\Omega$ of $\mathbb{R}^{\mathscr{S}}_{\geq}$ containing the positive orthant $\mathbb{R}^{\mathscr{S}}_{>}$. Subsets of power law kinetics though can be defined for proper supersets of $\mathbb{R}^{\mathscr{S}}_{>}$. We introduce a concept to describe this property:

\begin{definition}
The maximal definition domain $\Omega_{max}$ of a kinetics $K: \Omega \rightarrow \mathbb{R}^{\mathscr{R}}$ is the largest superset of $\Omega$ in $\mathbb{R}^{\mathscr{S}}_{\geq}$ on which $K$ is definable. A kinetics $K: \Omega \rightarrow \mathbb{R}^{\mathscr{R}}$ is positive orthant restricted if $\Omega_{max} = \mathbb{R}^{\mathscr{S}}_{>}.$
\end{definition}
Clearly, a power law kinetics is POR if and only if there is at least one negative value in each column of its kinetic order matrix. Non-POR power law kinetics include the subset of non-inhibitory kinetics (denoted by PL-NIK) with only nonnegative values in the kinetic order matrix. Other examples of non-POR  kinetics are mass action kinetics and Hill-type kinetics. For the latter, for each negative exponent $\mu$, one can write $\frac{x^\mu}{d + x^\mu} = \frac{1}{dx^{|\mu|} + 1}$, which then enables zero values for the species. 

\subsection{KSE Kinetics, a partial converse and a necessary condition}

We introduce the key property for our considerations:

\begin{definition}
A kinetic system $\left( \mathscr{N},K \right)$ has \textbf{kernel spanning equilibria images} (KSE kinetics) if $\langle K(E_+(\mathscr{N},K)) \rangle = \ker N$.
\end{definition}

Since $\langle K(E_+(\mathscr N,K))\rangle\subset \ker\mathscr N$, an equivalent condition is dim  $\langle K(E_+(\mathscr N,K))\rangle = r-s$, where $s$ is the $\dim S$ or $\dim (\Ima N)$. For POR kinetics, we have a simpler relationship:
 
 \begin{proposition}
 For a POR kinetics, $K(E_+\mathscr{N},K) = \ker N \cap \Ima K.$
 \end{proposition}
 \textbf{Proof. }
 As mentioned above, for any kinetics, $K(E_+(\mathscr{N},K)) \subset  \ker N \cap \Ima K \subset \ker N$. For $z \in \ker N \cap \Ima K$, $z = K(x)$, and since $K$ is POR, $x>0.$ Together with $N(z) = N(K(x)) = 0,$ this implies that $x \in E_+(\mathscr{N},K),$ and hence $z \in E_+(\mathscr{N},K).$

\begin{remark}\label{rs}
The above implies that for POR kinetics, we have $\langle K(E_+(\mathscr N,K))\rangle= \langle \ker N\cap \Ima K\rangle$. Hence to show that the kinetics is also KSE, we need only to construct a basis of $\ker N$ which is contained in $\Ima K.$  
\end{remark}

The partial converse is the following result:
\begin{theorem}\label{theorem:mainresult}
Let $\mathscr{N}$ be a weakly reversible CRN with a KSE kinetics $K$. If $(\mathscr{N},K)$ is absolutely complex balanced, then $\delta = 0$.
\end{theorem}
 
\textbf{Proof. }For any kinetic system, $K(Z_+(\mathscr{N},K)) \subset  \ker I_a \cap \Ima K \subset \ker I_a.$ Similarly,  $K(E_+(\mathscr{N},K)) \subset  \ker N \cap \Ima K \subset \ker N$. It follows that $\langle K(E_+(\mathscr{N},K)) \rangle \subset \ker N$. If $K$ has the KSE property, then equality holds. On the other hand, $(\mathscr{N},K)$ absolutely complex balanced implies that $K(Z_+(\mathscr{N},K)) = K(E_+(\mathscr{N},K))$, so that $\langle K(Z_+(\mathscr{N},K)) \rangle = \langle K(E_+(\mathscr{N},K)) \rangle = \ker N$. Combining this with
$\langle K(Z_+(\mathscr{N},K))\rangle$ $\subset \ker I_a \subset \ker N$ implies $\ker I_a = \ker N$. Since $\dim \ker N = \dim \ker I_a + \delta$, it follows that $\delta = 0$.

We obtain the following necessary condition for the occurrence of an ACB system with positive deficiency:

\begin{corollary}
If an ACB system has positive deficiency, then it is not a KSE system.
\end{corollary}

\begin{remark}
A more modular version of the necessary condition is the following:
\begin{enumerate}
    \item [i.] If a kinetic system is ACB, then $\dim \langle K(E_+(\mathscr{N},K)) \rangle\leq \dim \ker I_a = r - n + l.$
    \item [ii.] If a CRN has positive deficiency, then $\dim \ker I_a < \dim \ker N = r - s.$
\end{enumerate}
\end{remark}

Hence, an ACB system with positive deficiency has dim $\dim \langle K(E_+(\mathscr{N},K)) \rangle < r - s,$ i.e. is not KSE. To identify ACB systems with positive deficiency for a given set of kinetics, we only need to consider non-KSE systems.

\subsection{A technique for finding counterexamples to ACB}

Another reformulation of the partial converse is clearly: if a kinetic system has positive deficiency and is KSE, then it is not ACB. This enables the construction of non-ACB systems (if possible) for any given set of kinetics. We illustrate this with how we found the example in Section 3.

To show that this a non-ACB system, we use the KSE definition. Note that $K(x)$ is POR. By definition, we conclude the kinetics is KSE. Note that $r-s=5-1=4 = \dim \langle K(E_+(\mathscr{N},K)) \rangle =4$ ($K_2 = K_3$). 
  
 
 
 
 Using the contrapositive of Theorem \ref{theorem:mainresult}, since $\delta =2\neq 0$, then $(\mathscr{N},K)$ is not ACB. \\

\section{Summary}
\label{sec:summary} 
To conclude, we summarize our main results. 
\begin{enumerate}
    \item We revived the study on absolute complex balanced (ACB) systems by recalling results of Feinberg on zero deficiency systems (Feinberg ACB Theorem) and of Horn and Jackson on mass action systems (Horn and Jackson ACB Theorem). We tried to extend the latter problem by formulating necessary and sufficient conditions that ensure absolute complex balancing in systems with positive  deficiency (beyond those with mass action kinetics). 
    
    \item We showed an example of a complex balanced PL-RDK system which is not ACB leading to the non-extension of the ACB problem to PL-RDK systems.
    
    \item We established a necessary and sufficient condition for absolute complex balancing in CLP systems. 
    
    \item For non-CLP systems, we constructed and illustrated new classes of complex balanced systems that are absolutely complex balanced and with positive deficiency. The techniques used are known as the Positive Function Factor (PFF) method and the Coset Intersection Count (CIC) method. We also combine incidence independent and independent decompositions.
    
    \item We formulated a partial converse to the Feinberg ACB Theorem by identifying the set of KSE systems. Absolutely complex balanced KSE systems have zero deficiency. This enables us to construct a non-ACB system, that is, a KSE system with positive deficiency.
\end{enumerate}

\newpage
\appendix
\section{Notation and Acronym}
We denote by $\mathbb{R}$ and $\mathbb{Z}$ the set of real numbers and integers, respectively. For integers $a$ and $b$, let $\overline{a,b}= \{ j \in \mathbb{Z} | a \leq j \leq b \}.$ We denote the non-negative real numbers by $\mathbb{R}_{\geq0}$, and the positive real numbers by $\mathbb{R}_{>0}$. The sets $\mathbb{R}_{\geq 0}^n$ and $\mathbb{R}_{>0}^n$ are called the \textit{non-negative} and \textit{positive orthants} of $\mathbb{R}^n$, respectively. For $x \in \mathbb{R}^n$, the $i$th coordinate of $x$ is denoted by $x_i$, where $i \in \overline{1,n}$. We denote the vector space generated by vectors $y's$ as  $\langle y's \rangle$.

\begin{table}[h!]
\caption{List of symbols}
\begin{center}
\begin{tabular}{|l|c|}
\hline
Meaning & Symbol\\
\hline
augmented $T$ matrix & $\widehat T$\\
\hline
deficiency of a CRN & $\delta$\\
\hline
factor map of a kinetics $K$ & $\psi_K$ \\
\hline
incidence matrix of a CRN & $ I_a$\\
\hline 
$k-$ Laplacian matrix & $A_k$\\
\hline
kinetic flux subpace & $\widetilde{S}$\\
\hline
kinetic order matrix & $F$\\
\hline
kinetic order subspace & $\widetilde{S}$\\
\hline
kinetic reactant flux subspace & $\widetilde S_R$\\
\hline
kinetics of a CRN & $K$\\
\hline
kinetic reactant deficiency & $\widehat{\delta}$ \\
\hline
matrix of complexes  & $Y$\\
\hline
set of positive equilibria of  a system & $E_+(\mathscr N, K)$\\
\hline
set of complex balanced equilibria of  a system & $Z_+(\mathscr N, K)$\\
\hline
stoichiometric  matrix & $N$\\
\hline 
stoichiometric subspace of a CRN & $S$\\
\hline 
T matrix & $T$\\
\hline
\end{tabular}
\end{center}
\end{table}

\begin{table}[h!]
\caption{List of abbreviations}
\begin{center}
\begin{tabular}{|l|l|}
\hline
Abbreviation & Meaning\\
\hline
ACB & absolutely complex balanced \\
\hline
ACR & absolute concentration robustness \\
\hline
BCR & balanced concentration robustness \\
\hline
BIU & Bi-Independent Union \\
\hline
CRN & chemical reaction network\\
\hline
CIC & Coset Intersection Count \\
\hline
GMAK & generalized mass action kinetics\\
\hline
KSE & kernel spanning equilibria images \\
\hline
MAK & mass action kinetics \\
\hline
PFF & Positive Function Factor method \\
\hline
POR & positive orthant restricted \\
\hline
PLK & power-law kinetics\\
\hline
PL-FSK & power-law factor span surjective kinetics\\
\hline
PL-NIK & power law non-inhibitory kinetics\\
\hline
PL-RDK & power-law reactant determined kinetics\\
\hline
PL-TIK & $\widehat{T}$-rank maximal PL-RDK kinetics \\
\hline
PYK & Poly-PL kinetic systems \\
\hline
STAR-MSC & S-invariant Transformation by Adding Reactions \\
& via Maximal Stoichiometric Coefficients \\
\hline
SFRF & species formation rate function\\
\hline
\end{tabular}
\end{center}
\end{table}

\newpage
\section*{Acknowledgments}
ECJ acknowledges the funding by the UP System Enhanced Creative Work and Research Grant (ECWRG-2020-1-7-R) for this work.
\baselineskip=0.25in
\bibliographystyle{unsrt}  

\end{document}